\documentclass[11pt]{amsart}

\usepackage[latin1]{inputenc}

\usepackage{amssymb}
\usepackage{graphicx}
\usepackage{color}

\usepackage{hyperref}

\newcommand{\N}{\mathbb{N}}
\newcommand{\Z}{\mathbb{Z}}
\newcommand{\Q}{\mathbb{Q}}
\newcommand{\R}{\mathbb{R}}
\newcommand{\C}{\mathbb{C}}
\newcommand{\F}{\mathbb{F}}

\newcommand{\Aa}{\mathcal{A}}
\newcommand{\Bb}{\mathcal{B}}
\newcommand{\Cc}{\mathcal{C}}
\newcommand{\Dd}{\mathcal{D}}
\newcommand{\Ee}{\mathcal{E}}
\newcommand{\Ss}{\mathcal{S}}
\newcommand{\Ll}{\mathcal{L}}
\newcommand{\Mm}{\mathcal{M}}
\newcommand{\Nn}{\mathcal{N}}
\newcommand{\Xx}{\mathcal{X}}

\newcommand{\jj}{\mathcal{J}}
\newcommand{\Uu}{\mathcal{U}}

\newcommand{\om}{\omega}

\newcommand{\omuc}{\omega_{\mu,c}}

\newcommand{\J}{\mathcal{J}}

\newcommand{\Ssuc}{\Ss_{\mu,c}}
\newcommand{\Aauc}{\Aa_{\mu,c}}

\newcommand{\Gr}{\mathrm{Gr}}
\newcommand{\PD}{\mathrm{PD}}

\newcommand{\ev}{\mathrm{ev}}

\newcommand{\PU}{\mathrm{PU}}
\newcommand{\U}{\mathrm{U}}
\newcommand{\SO}{\mathrm{SO}}
\newcommand{\PSL}{\mathrm{PSL}}

\newcommand{\Symp}{\mathrm{Symp}}
\newcommand{\Diff}{\mathrm{Diff}}
\newcommand{\BSymp}{\mathrm{BSymp}}
\newcommand{\Aut}{\mathrm{Aut}}

\newcommand{\jjuc}{\jj_{\mu,c}}
\newcommand{\jjuco}{\jj^0_{\mu,c}}

\newcommand{\jjuci}{\jj^i_{\mu,c}}

\newcommand{\tGuc}{\widetilde{G}_{\mu,c}}

\newcommand{\Mdual}{M^{i,c}_{\mu}}

\newcommand{\Emb}{\mathrm{Emb}}
\newcommand{\IEmb}{\Im\mathrm{Emb}}

\newcommand{\CP}{\mathbb{C}P}
\newcommand{\tCP}{\widetilde{\CP}}
\newcommand{\STS}{S^2\times S^2}
\newcommand{\PbP}{\CP^2\#\,\overline{\CP}\,\!^2}
\newcommand{\NTB}{S^{2}\tilde{\times} S^{2}}
\newcommand{\tX}{\widetilde{X}}
\newcommand{\tM}{\widetilde{M}}
\newcommand{\Muo}{M^0_{\mu}}
\newcommand{\Mul}{M^1_{\mu}}
\newcommand{\Mui}{M^i_{\mu}}
\newcommand{\tMuco}{\widetilde{M}^0_{\mu,c}}
\newcommand{\tMucl}{\widetilde{M}^1_{\mu,c}}
\newcommand{\tMuci}{\widetilde{M}^i_{\mu,c}}

\newcommand{\into}{\hookrightarrow}
\newcommand{\lra}{\leftrightarrow}

\theoremstyle{plain}
\newtheorem{thm}{Theorem}[section]
\newtheorem*{thm*}{Theorem}
\newtheorem{prop}[thm]{Proposition}
\newtheorem*{prop*}{Proposition}
\newtheorem{lemma}[thm]{Lemma}
\newtheorem*{lemma*}{Lemma}
\newtheorem{cor}[thm]{Corollary}
\newtheorem*{cor*}{Corollary}

\theoremstyle{definition}
\newtheorem{defn}[thm]{Definition}
\newtheorem*{defn*}{Definition}
\newtheorem*{ackn}{Acknowledgments}

\theoremstyle{remark}
\newtheorem{remark}[thm]{Remark}
\newtheorem*{remark*}{Remark}

\newtheorem*{remarks*}{Remarks}

\begin{document}


\title[Symplectomorphism groups and embeddings]{Symplectomorphism groups and embeddings of balls into rational ruled $4$-manifolds}
\author{Martin Pinsonnault}
\email{mpinsonn@fields.utoronto.ca}
\address{Department of Mathematics, University of Toronto, Toronto, 
Canada, M5S 3G3.}
\curraddr{Fields Institute, Toronto, Canada, M5T 3J1.}
\thanks{This research was funded partly by NSERC grant BP-301203-2004.}


\begin{abstract}
Let $\Muo$ denote $S^{2}\times S^{2}$ endowed with a split symplectic form $\mu\sigma\oplus\sigma$ normalized so that $\mu\geq 1$ and $\sigma(S^{2})=1$. Given a symplectic embedding $\iota:B_{c}\into \Muo$ of the standard ball of capacity $c\in(0,1)$ into $\Muo$, consider the corresponding symplectic blow-up $\tMuco$. In this paper, we study the homotopy type of the symplectomorphism group $\Symp(\tMuco)$ and that of the space $\IEmb(B_{c},\Muo)$ of unparametrized symplectic embeddings of $B_{c}$ into $\Muo$. Writing $\ell$ for the largest integer strictly smaller than $\mu$, and $\lambda \in (0,1]$ for the difference $\mu-\ell$, we show that the symplectomorphism group of a blow-up of ``small'' capacity $c<\lambda$ is homotopically equivalent to the stabilizer of a point in $\Symp(\Muo)$, while that of a blow-up of ``large'' capacity $c\geq\lambda$ is homotopically equivalent to the stabilizer of a point in the symplectomorphism group of a nontrivial bundle $\PbP$ obtained by blowing down $\tMuco$. It follows that for $c<\lambda$, the space $\IEmb(B_{c},\Muo)$ is homotopy equivalent to $S^{2}\times S^{2}$, while for $c\geq \lambda$, it is not homotopy equivalent to any finite CW-complex. A similar result holds for symplectic ruled manifolds diffeomorphic to $\PbP$. By contrast, we show that the embedding spaces $\IEmb(B_{c},\CP^{2})$ and $\IEmb(B_{c_{1}}\sqcup B_{c_{2}},\CP^{2})$, if non empty, are always homotopy equivalent to the spaces of ordered configurations $F(\CP^{2},1)\simeq\CP^{2}$ and $F(\CP^{2},2)$. Our method relies on the theory of pseudo-holomorphic curves in $4$-manifolds, on the computation of Gromov invariants in rational $4$-manifolds, and on the inflation technique of Lalonde-McDuff.
\end{abstract}

\maketitle


\pagestyle{myheadings}
\markboth{SYMPLECTOMORPHISMS AND EMBEDDINGS}{MARTIN PINSONNAULT}

\section{Introduction}

Let $(X,\omega)$ be a symplectic $4$-manifold and consider the space $\Emb(B_{c},X)$ of symplectic embeddings of the standard ball $B_c\subset\R^4$ of radius $r$ and capacity $c=\pi r^2$ into $X$, endowed with its natural $C^{\infty}$ topology. We define the space of \emph{unparametrized} symplectic embeddings by setting
\[\IEmb(B_c, X):= \Emb(B_{c},X)/\Symp(B_{c})\]
The homotopy type of this space carries important information about the diffeomorphism class of the symplectic form $\om$. For instance, $\IEmb(B_c, X)$ is nonempty if and only if $c$ is smaller than the Gromov width of $(X,\om)$, and two symplectic blow-ups of $(X,\om)$ of same capacity $c$ are isotopic if $\IEmb(B_c, X)$ is connected. It is therefore natural to investigate the homotopy type of this space, at least for some amenable families of symplectic $4$-manifolds.

In a previous paper with F. Lalonde~\cite{LP-Duke}, we proposed a general framework to study the homotopy type of $\IEmb(B_c, X)$ based on the correspondence between embeddings of balls and symplectic blow-ups. We showed that when $\IEmb(B_c,X)$ is non-empty and connected, the natural action of the symplectomorphism group $\Symp(X,\om)$ on $\IEmb(B_c,X)$ is transitive and defines a fibration 
\begin{equation}\label{EqnMainFibration}
\Symp(\tX_{\iota},\Sigma)\to \Symp(X,\om)\to \IEmb(B_c,X)
\end{equation}
whose fiber over an embedding $\iota:B_{c}\into X$ is homotopy equivalent to the subgroup $\Symp(\tX_{\iota},\Sigma)$ of symplectomorphisms of the symplectic blow-up along $\iota$ sending the exceptional divisor $\Sigma$ to itself. We used this to investigate the space of symplectic balls in the product 
$$M_{\mu}^{0}:=(S^{2}\times S^{2},\mu\sigma\oplus\sigma)$$
where $\sigma$ is an area form such that $\sigma(S^{2})=1$ and where the  parameter $\mu$ lies in the interval $[1,2]$. By computing the rational homotopy invariants of the symplectomorphism group of the blow-up of $\Muo$ at a ball of capacity $c\in(0,1)$, and using the computation of the rational cohomology ring of $\Symp(\Muo)$ by Abreu-McDuff~\cite{AM}, we showed that the homotopy type of the space $\IEmb(B_{c},\Muo)$ depends on the capacity of the ball $B_{c}$ in an essential way:
\begin{thm}[(Lalonde-Pinsonnault~\cite{LP-Duke})] 
\begin{enumerate}
\item If $\mu=1$, then, for any $c\in(0,1)$, the embedding space $\IEmb(B_{c},M_{\mu}^{0})$ is homotopy equivalent to $S^{2}\times S^{2}$.
\item If $\mu\in(1,2]$ and $c\in(0,\mu-1)$, then $\IEmb(B_c, M_{\mu}^{0})$ is homotopy equivalent to $S^{2}\times S^{2}$.
\item If $\mu\in(1,2]$ and $c\in [\mu-1,1)$, then the rational homotopy groups of $\IEmb(B_{c},M_{\mu}^{0})$ vanish in all dimensions except in dimensions $2, 3$, and $4$, in which cases we have $\pi_2 =\Q^2$, $\pi_3 = \Q^3$, and $\pi_4 = \Q$. Moreover, this space is not homotopically equivalent to a finite CW-complex.
\end{enumerate}
\end{thm}

The main purpose of the present paper is the completion of these calculations for arbitrary rational ruled symplectic $4$-manifolds, that is, for all symplectic $S^2$-bundles over~$S^2$. As a byproduct, we also describe the homotopy type of $\IEmb(B_{c},\CP^{2})$ and that of the space $\IEmb(B_{c}\sqcup B_{c'},\CP^{2})$ of two disjoint balls of capacities $c$ and $c'$ in $\CP^{2}$.

\subsection*{General setting and main results}

\subsubsection*{Symplectic birational equivalences} In order to state our results, we first review some essential facts about rational ruled symplectic $4$-manifolds and their blow-ups.

To begin, recall that the classification theorem of Lalonde-McDuff~\cite{LM} asserts that any two cohomologous symplectic forms on a ruled $4$-manifold are symplectomorphic. It follows that a rational ruled symplectic $4$-manifold is, after rescaling, symplectomorphic to either
\begin{itemize}
\item the trivial bundle $\Muo := (S^2\times S^2, \om^0_{\mu})$, where the symplectic area of the a section $S^2\times\{*\}$ is $\mu\geq 1$ and the area of a fiber $\{*\}\times S^{2}$ is $1$; or
\item the non trivial bundle $\Mul := (\NTB,\om^{1}_{\mu})$, where the symplectic area of a section of self-intersection $-1$ is $\mu>0$ and the area of a fiber is $1$.
\end{itemize}
We will denote by $F^{i}$ the homology class of a fiber in $\Mui$, $i\in\{0,1\}$, and by $B^{i}$ the class of a section of self-intersection $-i$. With this choice of normalization, the area of a fiber is less or equal than the area of any embedded symplectic sphere of nonnegative self-intersection. It follows that the Gromov width of $\Mui$ is always $1=\om_{\mu}(F^{i})$, so that the embedding space $\IEmb(B_{c},\Mui)$ is nonempty if and only if $c\in(0,1)$. 

In~\cite{MD-Isotopy} McDuff showed that the symplectomorphism type of a symplectic blow-up of $\Mui$ along an embedded ball of capacity $c\in(0,1)$ depends only on the capacity $c$, and not on the particular embedding used in the construction (this is often referred to as ``the uniqueness of symplectic blow-ups'' for $\Mui$). Hereafter we will write $\tMuci$ for the symplectic blow-up of $\Mui$ at a ball of capacity~$c$, $\Sigma^{i}$ for the corresponding exceptional divisor, and $E^{i}=[\Sigma^{i}]$ for its homology class.

Now recall that the (smooth) blow-up of $M^0=S^2\times S^2$ can be identified with the blow-up of $M^1=\NTB$ via a diffeomorphism which preserves the fibers away from the exceptional loci. At the homology level, this identification induces an isomorphism $H_{2}(\widetilde{M}^{0};\Z)\simeq H_{2}(\widetilde{M}^{1};\Z)$ which acts on the classes $\{B^{0}, F^{0}, E^{0}\}$ and $\{B^{1},F^{1},E^{1}\}$ as follows: 
\begin{equation}\label{CorrespondanceHomologique}
\begin{array}{ccc}
F^{1}             & \lra & F^{0} \\
B^{1}+F^{1}-E^{1} & \lra & B^{0} \\
F^{1}-E^{1}       & \lra & E^{0} \\
B^{1}             & \lra & B^{0}-E^{0} \\
E^{1}             & \lra & F^{0}-E^{0}
\end{array}
\end{equation}
Note that $\widetilde{M}^{0}$ contains exactly $3$ exceptional classes (that is, classes represented by smooth embedded spheres of self-intersection $-1$), namely $E^{0}$, $F^{0}-E^{0}$, and $B^{0}-E^{0}$. The blow-down of an exceptional curve representing either $F^{0}-E^{0}$ or $B^{0}-E^{0}$ defines a manifold diffeomorphic to the non trivial bundle $\NTB$ while the blow-down of a curve in class $E^{0}$ yields $S^{2}\times S^{2}$. In particular, blowing up $S^{2}\times S^{2}$ and then blowing down along a curve representing $F^{0}-E^{0}$ defines a birational equivalence $\STS\dashrightarrow\NTB$ which preserves the fibers.

When one considers this birational equivalence in the symplectic category, the uniqueness of symplectic blow-ups implies that the blow-up of $\Muo$ at a ball of capacity $c\in(0,1)$ is symplectomorphic to the blow-up of $M^1_{\mu-c}$ at a ball of capacity $1-c$. Conversely, if $\mu\geq c$, the blow-up of $\Mul$ with capacity $c\in(0,1)$ is symplectomorphic to the blow-up of $M^0_{1+\mu-c}$ with capacity $1-c$. For this reason, given a ruled manifold $\Mui$ and a number $c\in (0,1)$ such that $c\leq\mu$, it is convenient to define the ``$c$-related'' manifold $\Mdual$ by setting
\begin{equation}\label{DefinitionDual}
\Mdual := 
\begin{cases}
M^1_{\mu-c} & \text{if $i=0$,}\\ 
M^0_{1+\mu-c} & \text{if $i=1$.} 
\end{cases}
\end{equation}
so that we can view the rational manifold $\tMuci$ either as a $c$-blow-up of $\Mui$ or as a $(1-c)$-blow-up of $\Mdual$. As we will see later, this leads to a duality between symplectic balls in $\Mui$ and in $\Mdual$. 

In the special case of $\Mul$ with $0<\mu<c<1$, the blow-up $\tMucl$ is again symplectomorphic to $\widetilde{M}_{1+\mu-c,1-c}^{0}$. However, the base class $B$ now has smaller area than the fibre class $F$, so that our choice of normalization forces us to interchange the two $S^{2}$ factors in the product $S^{2}\times S^{2}$ and to rescale the symplectic form. On the blow-up manifold, this induces a diffeomorphism which swaps the exceptional classes $B^{1}=B^{0}-E^{0}$ and $E^{1}=F^{0}-E^{0}$. It follows that $\tMucl$ is conformally symplectomorphic to $\widetilde{M}_{\mu',c'}^{0}$ where $\mu':=\frac{1}{1+\mu-c} \geq 1$ and $c':=\frac{1-c}{1+\mu-c}\in(0,1)$. Of course, this lack of symmetry will be reflected in the structure of $\IEmb(B_{c},\Mui)$ when $0<\mu<c<1$.

\subsubsection*{Symplectomorphism groups of $\tMuci$} The greater part of this paper is devoted to the determination of the homotopy type of the symplectomorphism groups $\Symp(\tMuci)$ and $\Symp(\tMuci,\Sigma)$. As in the works of Gromov~\cite{Gr}, Abreu~\cite{Ab}, and Abreu-McDuff~\cite{AM}, our understanding comes from a study of the natural action of these groups on the contractible space $\jjuci$ of tamed almost complex structures on $\tMuci$. This analysis was initiated in~\cite{LP-Duke} where we showed that this action preserves a stratification of $\jjuci$ whose strata are indexed by certain configurations of $J$-holomorphic curves or, equivalently, by equivalence classes of toric structures (that is, conjugacy classes of 2-tori $T^{2}\subset\Symp(\tMuci)$). For instance, if we consider the blow-up of $\Mui$ at a ball of capacity $c\leq\mu$, and if we write $\mu=\ell+\lambda$ with $\ell\in\N$ and $\ell<\mu\leq\ell+1$, the number $N$ of toric structures on $\tMuci$ is 
\[
N=
\begin{cases}
2\ell+1+i & \text{if $c<\lambda$,}\\ 
2\ell+i& \text{if $\lambda\leq c$.} 
\end{cases}
\]
so that the geometry of the stratification changes precisely when $\ell$ crosses an integer or when $c$ crosses the critical number $\lambda$. 

The analysis of $J$-holomorphic curves in $\tMuci$ is carried further in Section~\S\ref{SectionStructureCourbes} of the present paper. This is used in Section~\ref{SectionStabilite} to study how the homotopy types of the groups $\Symp(\tMuci)$ and $\Symp(\tMuci,\Sigma)$ change as the parameters $\mu$ and $c$ vary. We first note that in most cases, the distinction between $\Symp(\tMuci)$ and its subgroup $\Symp(\tMuci,\Sigma)$ is immaterial:

\begin{lemma}\label{RetractionSymp}
Given $c\leq\mu$, let $S$ be an embedded symplectic sphere in $\tMuci$ representing one of the exceptional classes $E^{i}$ or $F^{i}-E^{i}$. Then $\Symp(\tMuci)$ is homotopy equivalent to its subgroup $\Symp(\tMuci,S)$. 
\end{lemma}

Our main result states that, assuming $c\leq\mu$, the homotopy type of $\Symp(\tMuci,\Sigma)$ changes only when $\ell$ or $c$ cross a critical number, that is, when the set of toric stuctures on $\tMuci$ changes. This is the content of the following ``stability'' theorem:

\begin{thm}\label{StabilityTheorem}
Fix $i\in\{0,1\}$. Let $\mu_{1}\geq 1-i$ and write $\mu_{1}=\ell+\lambda_{1}$ with $\ell\in\N$ and $\lambda_{1}\in(0,1]$. Similarly, given a number $\mu_{2}\geq 1-i$ in the interval  $(\ell,\ell+1]$, write $\mu_{2}=\ell+\lambda_{2}$ with $\lambda_{2}\in(0,1]$. Consider $c_{1},c_{2}\in(0,1)$ such that either $c_{1}<\lambda_{1}$ and $c_{2}<\lambda_{2}$, or $\lambda_{1}\leq c_{1}$ and $\lambda_{2}\leq c_{2}$. Suppose also that $c_{i}\leq\mu_{i}$. Then the symplectomorphism groups $\Symp(\widetilde{M}^i_{\mu_{1},c_{1}},\Sigma)$ and $\Symp(\widetilde{M}^i_{\mu_{2},c_{2}}, \Sigma)$ are homotopy equivalent.
\end{thm}

Observe that from the homotopy fibration~(\ref{EqnMainFibration}) above, it follows that the behaviour of symplectic balls of capacity $c\leq\mu$ in $\Mui$ depends only on whether $c<\lambda$ or $\lambda\leq c$. This motivates the following definition: 

\begin{defn}\label{DefinitionSmall-Big}
Write $\mu=\ell+\lambda$ with $\ell\in\N$ and $\ell<\mu\leq\ell+1$. We say that an embedded symplectic ball $B_{c}\into \Mui$ is {\em small} if its capacity $c$ is strictly less than $\lambda$. It is said to be {\em big} if $\lambda\leq c$.
\end{defn}

It is easy to see that if $c\leq\mu$, the blow-up of a ``small'' ball $B_c$ in $\Mui$ is  symplectomorphic to the blow-up of a ``big'' ball in the $c$-related manifold $\Mdual$. Conversely, blowing up a ``big'' ball in $\Mui$ is equivalent to blowing up a ``small'' ball in $\Mdual$. Assuming $c<\lambda$ and looking at the limit $c\to 0$ (while keeping $\mu$ constant), we obtain the following simple homotopy-theoretic description of the groups $\Symp(\tMuci,\Sigma)\simeq\Symp(\tMuci)$\,: 

\begin{thm}\label{HomotopyOfSymp} 
Let $\tMuci$ be the symplectic blow-up of $\Mui$ at a ball of capacity $c\in(0,1)$. In the case $i=1$, suppose also that $c\leq\mu$. Let write $\mu=\ell+\lambda$ where $\ell$ is the unique integer such that $\ell<\mu\leq\ell+1$ and where $\lambda\in(0,1]$. Then,
\begin{enumerate}
\item when $c<\lambda$, the symplectomorphism group of the blow-up manifold $\tMuci$ is homotopy equivalent to the group $\Symp(\Mui,p)$ of symplectomorphims that fix a point $p$ in $\Mui$. 
\item When $\lambda\leq c$, the symplectomorphism group of the blow-up manifold $\tMuci$ is homotopy equivalent to the group $\Symp(\Mdual,p)$ of symplectomorphims that fix a point $p$ in the $c$-related manifold $\Mdual$. 
\end{enumerate}
\end{thm}

In the case $0<\mu<c<1$, the symplectomorphism group $\Symp(\tMucl)$ is homeomorphic to the group $\Symp(\widetilde{M}_{\mu',c'}^{0})$ where $\mu':=\frac{c}{1+\mu-c}$ and $c':=\frac{\mu}{1+\mu-c}$. However, the groups $\Symp(\tMucl)$ and $\Symp(\tMucl,\Sigma)$ are no longer homotopy equivalent. Nevertheless, the techniques used in the proof of Theorem~\ref{StabilityTheorem} apply mutatis mutandis and yield the following simple result:

\begin{thm}\label{HomotopyOfSympSpecialCase}
When $0<\mu<c<1$, the symplectomorphism group $\Symp(\tMucl,\Sigma)$ is homotopy equivalent to the stabilizer of a point in $\Mul$.
\end{thm}

From the homotopy equivalences of Theorem~\ref{HomotopyOfSymp} and Theorem~\ref{HomotopyOfSympSpecialCase}, and from the computations of the rational homotopy invariants of $\Symp(\Mui)$ by Abreu-McDuff~\cite{AM}, we easily obtain a description of the rational homotopy module $\pi_{*}(\Symp(\tMuci))\otimes\Q$, of the cohomology ring $H^{*}(\Symp(\tMuci);\Q)$, of the module $H^{*}(\BSymp(\tMuci);\Q)$, and of the Pontrjagin and Samelson products in $H_{*}(\Symp(\tMuci);\Q)$ and $\pi_{*}(\BSymp(\tMuci))\otimes\Q$, see~\S\ref{SectionFormalConsequences}.  

\subsubsection*{Spaces of symplectic embeddings} Our understanding of the symplectomorphism groups of $\tMuci$ can now be applied to the study of the embedding space $\IEmb(B_{c},\Mui)$ via the homotopy fibration~(\ref{EqnMainFibration}). In the case $c\leq\mu$, this gives:

\begin{thm}\label{HomotopyOfEmbeddings}
Given a rational ruled manifold $\Mui$, let write $\mu=\ell+\lambda$ where $\ell$ is the unique integer such that $\ell<\mu\leq\ell+1$ and where $\lambda\in(0,1]$. Let $c\in(0,1)$. In the case $i=1$, suppose also that $c\leq\mu$. Then,
\begin{itemize}
\item when $c<\lambda$, the space $\IEmb(B_c,\Mui)$ is homotopy equivalent to $\Mui$, that is, small symplectic balls in $\Mui$ behave like points.
\item When $\lambda\leq c$, the rational homotopy of the space $\IEmb(B_c, \Mui)$ is generated, additively, by two elements of degree $2$, two elements of degree $3$, one element of degree $4\ell-1+2i$, and one element of degree $4\ell+2i$.
\end{itemize}
\end{thm}

In the case $0<\mu<c<1$, we obtain an even simpler result:

\begin{thm}\label{HomotopyOfEmbeddingsSpecialCase}
If $0<\mu<c<1$, the space $\IEmb(B_c,\Mul)$ is homotopy equivalent to $\Mul$.
\end{thm}

\begin{cor}
For all $\mu\in(0,1]$ and $c\in(0,1)$, $\IEmb(B_c,\Mul)$ is homotopy equivalent to $\Mul$.
\end{cor}

The same techniques can be applied to investigate the homotopy theoretic properties of other embedding spaces. In particular, our computations allow us to describe the space $\IEmb(B_{c};\CP^{2})$ and the space of symplectic embeddings $B_{c_{1}}\sqcup B_{c_{2}}\into\CP^{2}$ of two disjoint balls in~$\CP^{2}$\,:

\begin{thm}\label{EmbeddingsInCP2}
Consider $\CP^{2}$ endowed with its standard Fubini-Study symplectic form normalized so that the area of a line is $1$. Then, 
\begin{enumerate}
\item for any $c\in(0,1)$, the space $\IEmb(B_{c},\CP^{2})$ is homotopy equivalent to $\CP^{2}$.
\item Given $c_{1}, c_{2}\in(0,1)$ so that $c_{1}+c_{2}<1$, the space $\IEmb(B_{c_{1}}\sqcup B_{c_{2}};\CP^{2})$ is homotopy equivalent to the space $F(2,\CP^{2})$ of ordered configurations of $2$ points in $\CP^{2}$.
\end{enumerate}
\end{thm}

We conclude this Introduction with some comments on possible generalizations of our results. Following Anjos-Granja~\cite{AG} and Abreu-Granja-Kitchloo~\cite{AGK}, the action of the symplectomorphism group $\Symp(\tMuci)$ on the natural stratification of the space of almost complex structures $\jjuci$ shall lead to a homotopy decomposition of $\BSymp(\tMuci)$ as a homotopy colimit over the poset category of toric structures and their intersections. In turn, this should yield a better description of the homotopy type of the embedding space $\IEmb(B_c,\Mui)$ in the case of ``big'' balls. In another direction, we note that if we impose certain conditions on the capacities, it is possible to get a partial understanding of the space $\IEmb(B_{c_{1}}\sqcup\cdots\sqcup B_{c_{r}},\CP^{2})$ of $r\geq 3$ disjoint symplectic balls in $\CP^{2}$ by proving versions of the Stability Theorem~\ref{StabilityTheorem}, see~\cite{Pi-EmbeddingsCP2}. In these manifolds, the structure of pseudo-holomorphic curves can be quite complicated and is best understood through the action of the group of birational equivalences of $X_{r}$. But new phenomena also occur which make the analysis of the symplectomorphism groups more difficult. For instance, a proper understanding of the action of $\Symp(X_{r})$ on the contractible space $\jj(X_{r})$ requires a description of various moduli spaces of complex structures, and this may no longer yield a computable description of the homotopy type of $\Symp(X_{r})$. Another problem is that the symplectomorphism groups $\Symp(X_{c_{1},\ldots,c_{r}})$ are not necessarily connected. Indeed, Seidel showed that for the \emph{monotone} blow-ups $X_{r}$, $5\leq r\leq 8$, the group $G=\Symp(X_{r})\cap\Diff_{0}$ is not connected, see~\cite{Se}. His proof relies on the action of $\pi_{0}(G)$ on the Donaldson quantum category of $X_{r}$ and yields computable invariants only in the monotone case. However, it seems possible to combine our method with that of Seidel in order to get new insights on the homotopy type of $\Symp(X_{r})$ in the non-monotone case. These ideas will be pursued elsewhere.

\subsubsection*{Conventions}
We endow all diffeomorphism groups with the $C^{\infty}$ topology. A theorem of Milnor~\cite{Mi} implies that these topological groups are homotopy equivalent to finite or countable CW-complexes. Consequently, we do not distinguish between weak homotopy equivalences and genuine homotopy equivalences. Note also that our notation differs slightly from the one used in~\cite{AM} where the manifold $M^{i}_{\lambda}$ corresponds to our~$\Mui$ with $\mu=1+\lambda$. Finally, since the present paper is a continuation of~\cite{LP-Duke}, we refer the reader to that paper for the details of several arguments and proofs. 

\section{Structure of $J$-holomorphic curves in $\tMuci$}\label{SectionStructureCourbes}

In this section, we describe the structure of $J$-holomorphic spheres representing special homology classes in the blow-up $\tMuci$. This analysis relies on well-known results about $J$-holomorphic curves in symplectic $4$-manifolds that we briefly recall for convenience. The proofs can be found in~\cite{AL, Gr, HLS, MS}. 

\subsection{$J$-holomorphic curves in symplectic $4$-manifolds}

A \emph{parametrized $J$-sphere} is a map $u:\CP^1 \to M$ which satisfies the Cauchy-Riemann equation $du \circ j = J \circ du$. We always assume that $u$ is \emph{somewhere injective}, that is, there is a point $z\in \CP^1$ such that $du_{z}\neq 0$ and $z=u^{-1}(u(z))$. In this case, $u$ is not a multiple cover and its image $S=u(\CP^{2})$ represents a homology class in $H_{2}(X;\Z)$ that we denote $[S]$ or $[u]$.

An almost complex structure, $J$, is \emph{tamed} by a symplectic form $\omega$
if $\omega(v,Jv) > 0$ for all $v \neq 0$. Note that if $J \in \J$, then every $J$-holomorphic submanifold $C \subset M$ is symplectic and, conversely, every symplectic submanifold of $(M,\om)$ is $J$-holomorphic for some $J\in\J$. In general, if a class $A\neq 0$ can be represented by a $J$-holomorphic sphere, for some $J\in\J(\om)$, then $\om(A) > 0$. For any symplectic manifold, the space of tamed almost complex structures, denoted by $\J=\J(M,\omega)$, is always non-empty and contractible. It follows that the first Chern class $c_1(TM)$ is independent of the choice of $J \in \J$.

Given a class $A\in H_{2}(X;\Z)$, we define the universal moduli space $\Mm_{p}(A,\J)$ of parametrized holomorphic $A$-spheres by setting:
$$\Mm_{p}(A,\J) = \{ (u,J)~~:~~J\in\J, \text{~$u$ is $J$-holomorphic, and~} [u]=A \}$$
A fundamental result of the theory of $J$-holomorphic curves asserts that $\Mm_{p}(A,\J)$ is always a Fréchet manifold and that the projection $P_{A}:\Mm_{p}(A,\J)\to\J$ is a Fredholm map of index $2c_{1}(A)+4$. Given $J\in\J$, the preimage $\Mm_{p}(A,J)=P_{A}^{-1}(J)$ is the space of all parametrized $J$-holomorphic $A$-spheres. The automorphism group $G=\PSL(2,\C)$ of $\CP^1$ acts freely on $\Mm_{p}(A,J)$ by reparametrization and the quotient $\Mm(A,J)=\Mm_{p}(A,J)/G$ is the moduli space of un\-pa\-ra\-me\-tri\-zed $J$-holomorphic $A$-spheres. An almost complex structure, $J\in\J$, is said to be regular for the class $A$ if it is a regular value for the projection $P_A$. In this case, the sets $\Mm_p(A,J)$ and $\Mm(A,J)$ are smooth manifolds of dimensions $2c_1(A)+4$ and $2c_{1}(A)-2$ respectively. The set $J_{\mathrm{reg}}$ of all such regular $J$ is a subset of second category in $\J$.

Although the moduli space $\Mm_{p}(A,J)$ is never compact, its quotient $\Mm(A,J)=\Mm_{p}(A,J)/G$ can be compactified by adding \emph{cusp-curves}, that is, connected unions of possibly multiply-covered $J$-curves\,\footnote{Although the compactification of the moduli space is best described using the finer notions of \emph{stable maps} and \emph{stable curves}, the simpler notion of cusp-curves is sufficient for our purpose.}. Indeed we have the following theorem:

\begin{thm*}[Gromov's Compactness Theorem]\label{Compacite} If a sequence $\{J_i\}\subset\J$ converges to $J_\infty\in\J$ and if $S_i$ are unparametrized $J_i$-holomorphic spheres of bounded symplectic area, then there is a subsequence of the $S_i$ which converges weakly to a genuine $J$-sphere or to a $J$-holomorphic cusp-curve.
\end{thm*}

A crucial fact is that weak convergence implies convergence in the Hausdorff topology of subsets. Consequently, the area of the limit is the limit of the areas of the converging subsequence. In particular, if all the $S_i$ represent the same class $A$, and if the sequence $S_{i}$ converge to a cusp-curve $S_\infty$ made of several components $C_1\cup\cdots\cup C_k$, each having multiplicity $m_{j}>0$, then their homology classes $m_{j}[C_j]$ provide a decomposition of $A$ with $\om(A)>\om(C_j)>0$. 

In dimension $4$, $J$-holomorphic spheres have a number of special properties which often reduce geometric questions to homological ones. These are:

\begin{trivlist}
\item {\sl Regularity:}\label{Regularity} 
The Hofer-Lizan-Sikorav regularity criterion asserts that if $c_{1}(A)\geq 1$, then all $J\in\J$ are regular values of the projection $P_{A}$.
\smallskip
\item {\sl Positivity of intersections:}\label{PositivityOfIntersections}
Two distinct $J$-spheres $S_1$ and $S_2$ in a $4$-manifold have only a finite number of intersection points. Each such point contributes positively to the algebraic intersection number $S_1\cdot S_2$. Moreover, $S_1\cdot S_2=0$ if and only if the spheres are disjoint, and $S_1\cdot S_2=1$ if and only if they meet exactly once transversally.
\smallskip
\item {\sl Adjunction formula:}\label{Adjunction}
Given a symplectic $4$-manifold $M$ and a class $A\in H_{2}(M;\Z)$, we define the virtual genus of $A$ by the formula
$$g_{v}(A)=1+\frac{1}{2}\left(A\cdot A - c_{1}(A)\right)$$
If $A$ is represented by the image of a somewhere injective $J$-holomorphic map $u:\CP^{1}\to M$, then $g_v(A)$ is a non-negative integer which is zero if and only if $u$ is an embedding.
\smallskip
\item {\sl Exceptional curves:}\label{ExceptionalCurves}
[see Lemma 3.1 in~\cite{MD-Structure}] 
Let $A \in H_2(M;\Z)$ be a homology class which is represented by an embedded symplectic sphere $S$, and such that $A \cdot A = -1$. Then the space $\J_{A}$ of all $J\in\J$ for which the class $A$ is represented by an embedded $J$-sphere contains an open, dense, and path-connected subspace of $\J$.
\end{trivlist}

\subsection{Gromov invariants} We now recall some facts from Taubes-Seiberg-Witten theory and Gromov invariants. Given a 4-dimensional symplectic manifold $(X,\om)$ and a class $A\in H_{2}(X;\Z)$ such that $k(A)=\frac{1}{2}(A\cdot A + c_{1}(A))\geq 0$, the Gromov invariant of $A$ as defined by Taubes\footnote{The original definition given by Taubes must be modified slightly to handle multiply covered curves and to have an equivalence between Gromov invariants and Seiberg-Witten invariants. See~\cite{LL:EquivalenceGr-SW} for a complete discussion.} in~\cite{Ta} counts, for a generic almost complex structure $J$ tamed by $\om$, the algebraic number of embedded $J$-holomorphic curves in class $A$ passing through $k(A)$ generic points. In general, these curves may be disconnected and, in that case, components of nonnegative self-intersection are embedded holomorphic curve of some genus $g\geq 0$, while all other components are exceptional spheres. An important corollary of Gromov's compactness theorem is that, given a tamed $J$, any class $A$ with non zero Gromov invariant $\Gr_{\om}(A)$ is represented by a collection of $J$-holomorphic curves or cusp-curves. 

The Gromov invariants only depend on the deformation class of $\om$, that is, $\Gr_{\om_{t}}(A)$ is constant along any $1$-parameter family of symplectic forms $\om_{t}$. In fact, in dimension four, it follows from the work of Taubes that Gromov invariants are also smooth invariants of $M$. 

Now suppose that $(X,\om)$ is some blow-up of a rational ruled manifold $\Mui$ with canonical homology class $K=-\PD(c_{1})$. Since $b_{2}^{+}=1$, the wall-crossing formula of Li-Liu for Seiberg-Witten invariants (see~\cite{LL}, Corollary 1.4) implies that for every class $A\in H_{2}(X,\Z)$ verifying the condition $k(A)=\frac{1}{2}(A\cdot A + c_{1}(A))\geq 0$ we have
\begin{equation}\label{SymetrieGr}
\Gr(A)\pm \Gr(K-A) = \pm 1
\end{equation}
In particular, if $\Gr(K-A)=0$ (for example, when $\om(K-A)\leq 0$), then $\Gr(A)\neq 0$. 

\subsection{Structure of $J$-holomorphic curves in $\tMuci$}

As before, we write $\mu=\ell + \lambda$ with $\ell\in\N$ and $\lambda\in(0,1]$. Since we can always obtain $\tMuci$ by blowing up a ball in the product $S^{2}\times S^{2}$, we will work only with $\Muo$ and, to simplify notation, we will omit the superscript ``$0$'' on homology classes in $H_{2}(\tMuco,\Z)$. In particular, we will write $B$ for $[S^{2}\times\{*\}]$, $F$ for $[\{*\}\times S^{2}]$, and $E$ for the the class of the exceptional divisor in $\tMuco$.

\begin{lemma}\label{PositivityOfp}
Given $\mu\geq 1$, $c\in(0,1)$, if a class $A=pB + qF - rE\in H_{2}(\Muo,\Z)$ has a simple $J$-holomorphic representative, then $p\geq 0$. Moreover, if $p=0$, then $A$ is either $F$, $F-E$ or $E$. Furthermore, if $p=1$, then $r\in\{0,1\}$.
\end{lemma}
\begin{proof}
This follows easily from the adjunction inequality $0\leq 2g_v(A)=2(p-1)(q-1)-r(r-1)$, the positivity of the symplectic area $\om(A)=\mu p+q-cr$, and our choice of normalization: $\mu=\om(B)\geq\om(F)=1>\om(E)=c>0$.
\end{proof}

\begin{lemma}\label{PersistanceEetF-E}
Given $\mu\geq 1$, $c\in(0,1)$ and any tamed almost complex structure $J\in\jjuc$, each of the exceptional classes $E$ and $F-E$ are represented by a unique embedded $J$-holomorphic curves. Consequently, if a class $A = pB + qF - rE$ has a simple $J$-holomorphic representative, then $p\geq r\geq 0$, unless $A\in\{E, F-E\}$.
\end{lemma}
\begin{proof}
We know that the set of almost complex structures $J$ for which the exceptional classes have embedded $J$-holomorphic representatives is open and dense. Moreover, by positivity of intersections, any such $J$-holomorphic representative is necessarily unique. Thus, by Gromov's compactness theorem, we only need to show that these classes cannot be represented by multiply covered curves nor by cusp-curves.

First, consider the class $E$. Since $E$ is a primitive class, it cannot be represented by a multiply covered curve. Suppose it is represented by a cusp curve $C=\bigcup_1^N m_iC_i$, where $N\geq 2$ and each component $C_i$ is simple. By lemma~\ref{PositivityOfp}, the homology class of each component $C_i$ is either $F$, $E$ or $F-E$. The cases $[C_i]=F$ and $[C_i]=E$ are impossible since the symplectic area of each component must be strictly less than the area of $E$. But then each component $C_i$ represents the class $F-E$, which is also impossible. Therefore, the class $E$ cannot be represented by a cusp-curve.

The same arguments apply to the class $F-E$ and show that it is always represented by an embedded $J$-holomorphic sphere. Consequently, if a class $A = pB + qF - rE$ has a simple $J$-holomorphic representative, the positivity of intersections implies that $0\leq A\cdot (F-E) = p-r$ and $0\leq A\cdot E=r$, unless $A$ is either $E$ or $F-E$.
\end{proof}

\begin{lemma}\label{ModuliSpaceFiber}
Given $\mu\geq 1$, $c\in(0,1)$ and a tamed almost complex structure $J\in\jjuc$, the moduli space $\Mm(F,J)$ of embedded $J$-holomorphic spheres representing the fiber class $F$ is a non-empty, open, $2$-dimensional manifold.
\end{lemma}
\begin{proof}
Since the class $F$ is regular $(c_1(F)=2\geq 1)$, the dimension of the space of unparametrized $J$-holomorphic representatives is $2c_1(F)-2=2$ whenever it is nonempty. But since $Gr(F)=1$, this space is non-empty for generic tamed almost complex structure $J$, and there is exactly one embedded $J$-holomorphic representative of $F$ passing through $k(F)=\frac{1}{2}(c_1(F) + F^2) = 1$ generic point of~$\tMuco$.

Now suppose that for some tamed $J$, there is no embedded $J$-holomorphic curve representing the class $F$. By Gromov's compactness theorem, a generic point of $\tMuco$ must belong to either a multiply-covered curve or a cusp-curve in class $F$. Being simple, the class $F$ cannot be represented by a multiply covered curve. If it is represented by a $J$-holomorphic cusp-curve, the two previous lemmas imply that this cusp-curve is necessarily of the form $(F-E)\cup E$ so that it cannot contain a generic point. We conclude that $\Mm(F,J)$ is nonempty for any tamed almost complex structure~$J$.
\end{proof}

Let us define the classes $D_{i}\in H_{2}(\tMuco,\Z)$, $i\in\Z$, by setting $D_{2k+1}=B+kF$ and $D_{2k}=B+kF-E$. Note that we have $D_{i}\cdot D_{j}=\lfloor\frac{i+j}{2}\rfloor$, where $\lfloor x \rfloor$ denotes the greatest integer strictly less than~$x$. In particular, $D_{i}\cdot D_{i}=i-1$ and the adjunction formula implies that any $J$-holomorphic sphere in class $D_{i}$ must be embedded. Note also that the virtual dimension of the spaces of unparametrized $J$-curves representing $D_{i}$ is $2i$ and that $k(D_{i})=i$. 

\begin{lemma}\label{Gr(Di)}
Given an integer $i\geq 1$, the Gromov invariant of the class $D_{i}$ is nonzero.
Consequently, for a generic tamed almost complex structure $J$, there exist an embedded $J$-holomorphic sphere representing $D_{i}$ and passing through a generic set of $k(D_{i})=i$ points in $\tMuco$.
\end{lemma}
\begin{proof}
Let $K=-2(B+F)-E$ be the Poincar\'e dual of the canonical class of $\tMuco$. Because $\om(K-D_{i})\leq 0$, we know that $\Gr(K-D_{i})=0$, and since $k=k(D_{i})\geq 0$, formula~(\ref{SymetrieGr}) implies that $\Gr(D_{i})=\pm 1$. Given a generic $J$, the class $D_{i}$ is thus represented by the disjoint union of some embedded $J$-holomorphic curve of genus $g\geq 0$ with some exceptional $J$-spheres. But since $D_{i}$ has nonnegative intersection with the three exceptional classes $B-E$, $F-E$, and $E$, this representative is in fact connected and its genus is given by the adjunction formula.
\end{proof}

\begin{lemma}\label{Stratification}
The set of tamed almost complex structures on $\tMuco$ for which the class $D_{-1}=B-E$ is represented by an embedded $J$-holomorphic sphere is open and dense in $\jjuc$. If for a given $J$ there is no such sphere, then there is a unique integer $1\leq m \leq \ell$ such that either the class $D_{-2m}=B-mF$ or the class $D_{-2m-1}=B-mF-E$ is represented by a unique embedded $J$-holomorphic sphere.
\end{lemma}
\begin{proof}
Since the class $B-E$ is exceptional, the set of almost complex structures $J$ for which it has an embedded $J$-holomorphic representative is open and dense. Being primitive, it has no multiply-covered representative. Suppose it is represented by a cusp curve $C=\bigcup_1^N m_iC_i$, where $N\geq 2$. We must show that one of its component is a representative of $B-mF$ or $B-mF-E$ for some integer $1\leq m \leq \ell$. By lemma~\ref{PositivityOfp} and lemma~\ref{PersistanceEetF-E}, any cusp-curve gives a decomposition of the class $B-E$ of the form 
\[
D_{i} = (B+mF-rE)  + \sum_a s_aF +  t(F-E) + uE
\]
where the multiplicities $s_a$, $t$, and $u$ are all positive and where
$r\in\{0,1\}$. This implies that $m<0$. Since the symplectic area of the component in class $B+mF-rE$ must be strictly positive, $|m|$ must be less or equal than $\ell$. Note that if $|m|<\ell$, the class $B-|m|F-E$ has positive area since $c<1$, while if $|m|=\ell$, it has positive area only if $c < \lambda$. Finally, since the intersection of $B-|m|F-|r|E$ with any other class $D_{i}$, $i\leq 0$, is negative, at most one such class can be represented by a $J$-holomorphic curve for a given tamed almost complex structure $J$.
\end{proof}

Let $\tilde{\jj}_{\mu,c,m}$ be the space of all almost complex structure $J\in\tilde{\jj}_{\mu,c}$ such that the class $D_{-m}$ is represented by an embedded $J$-holomorphic sphere. Lemma~\ref{Stratification} shows that the space $\tilde{\jj}_{\mu,c}$ is the disjoint union of its subspaces $\tilde{\jj}_{\mu,c,m}$, that is,

\begin{equation}\label{LaStratification}
\tilde{\jj}_{\mu,c}=\tilde{\jj}_{\mu,c,0}\sqcup\cdots\sqcup\tilde{\jj}_{\mu,c,N}
\text{~~~where~} N=
\begin{cases}
2\ell & \text{if $c < \lambda$,}\\
2\ell-1   & \text{if $\lambda\leq c$}
\end{cases}
\end{equation}

\begin{lemma}\label{ExistenceCourbesDi}
Given a tamed almost complex structure $J$ belonging to the strata $\tilde{\jj}_{\mu,c,m}^{0}$, a class $D_{i}$, $i\neq -m$, is represented by some embedded $J$-holomorphic sphere if and only if $D_{i}\cdot D_{-m}\geq 0$. In particular, if either $c<\lambda$ and $i\geq 2\ell+1$, or $\lambda\leq c$ and $i\geq 2\ell$, then, for all $J\in\tilde{\jj}_{\mu,c}$, the class $D_{i}$ has embedded $J$-holomorphic representatives.
\end{lemma}
\begin{proof}

We first note that the presence of a $J$-holomorphic curve in class $D_{-m}$ forbids the existence of an embedded representative of any class $D_{i}$, $i\neq -m$, such that $D_{i}\cdot D_{-m}<0$. Now given any $J\in\tilde{\jj}_{\mu,c,m}$, consider a class $D_{i}$ such that $D_{i}\cdot D_{-m}\geq 0$ and suppose that $D_{i}$ does not have any embedded $J$-holomorphic representative. Note that we have necessarily $i\geq 1$ and hence $\Gr(D_{i})\neq 0$. Since $D_{i}$ is generically represented by embedded holomorphic spheres, there must be a $J$-holomorphic cusp-curve representing $D_{i}$ and passing through $k=k(D_{i})=i$ generic points $\{p_{1},\ldots,p_{k}\}$. Again, Lemma~\ref{PositivityOfp} and lemma~\ref{PersistanceEetF-E} imply that any such $J$-holomorphic cusp-curve $C=\bigcup_1^N m_iC_i$ with at least two components must defines a decomposition of $D_{i}$ of the form
$$D_{i} = D_{r}  + \sum_a s_aF + \sum_b t_b(F-E) + \sum_c u_cE$$
where the multiplicities $s_a,t_b, u_c$ are all positive and where $D_{r}\cdot D_{-m}\geq 0$ unless $r=-m$.

Consider the reduced cusp-curve $\bar{C}$ obtained from $C$ by removing all but one copy of its repeated components, and by replacing multiply-covered components by their underlying simple curve. The curves $\bar{C}$ and $C$ have the same image and hence contain the same subset of special points in $\{p_{1},\ldots,p_{k}\}$. The cusp-curve $\bar{C}$ represents a class of the form
$$[\bar{C}] = D_{r} + \sum_a F + \sum_b (F-E) + \sum_c E$$
since, by the adjunction formula, there are no simple curves representing $mF$, $m(F-E)$, or $mE$ for $m\geq 2$. When $r\geq 0$, all the components are regular and the dimension of the space of unparametrized reduced cusp-curves of this type is simply
$$
\dim\,\Mm(\bar{C},J)=\sum_{\text{components $\bar{C}_{i}$}}\dim\,\Mm([\bar{C}_{i}],J)=2c_1(\bar{C})-2N
$$
where $N$ is the number of components of $\bar{C}$. Since $c_1(\bar{C})\leq c_1(D_{i})$, the dimension of unparametrized and reduced cusp-curves of type $\bar{C}$ is less than $2k(D_{i})$ and this implies that no such cusp-curve can pass through $k(D_{i})$ generic points.

When $r<0$, the only possibility is $r=-m$. The component in class $D_{-m}$ has negative self-intersection and is thus isolated. Therefore, the dimension of the space of cusp-curves of type $\bar{C}$ is at most $2n$, where $n$ is the number of components representing the class $F$. Note that this number $n$ is maximal when the reduced cusp-curve is of one of the three following forms:
$$
\begin{array}{ccc}
D_{i}=D_{-m}+\sum F, & D_{i}=D_{-m}+\sum F + (F-E), & \text{or~}D_{i}=D_{-m}+\sum F+E
\end{array}
$$ 
In each case the dimension of the moduli space of cusp-curves is strictly less than $2k(D_{i})=2i$ and, again, this implies that no cusp-curve of type $\bar{C}$ can pass through $k(D_{i})$ generic points.

The last statement of the lemma follows from the fact that when $c<\lambda$, the highest codimensional strata in $\tilde{\jj}_{\mu,c}$ corresponds to structures for which the class $D_{-2\ell}=B-\ell F-E$ is represented, and $D_{-2\ell}\cdot D_{i}\geq 0$ if and only if $i\geq 2\ell+1$. Similarly, when $\lambda\leq c$, the last strata corresponds to curves $D_{-2\ell+1}=B-\ell F$ and $D_{-2\ell+1}\cdot D_{i}\geq 0$ if and only if $i\geq 2\ell$.
\end{proof}

\begin{remark}\label{RemarkSpecialCase}
Our identification of $\tilde{M}_{\mu+1-c,1-c}^{0}$ with $\tilde{M}_{\mu,c}^{1}$ implies that when we blow up $M_{\mu}^{1}$ with capacity $0<\mu<c<1$, the exceptional classes that are always represented by embedded $J$-spheres are $B^{1}$ and $F^{1}-E^{1}$. In particular, the class $E^{1}$ of the exceptional divisor can be represented by $J$-holomorphic cusp-curves. This will force us to study the space $\IEmb(B_{c},M_{\mu}^{1})$, $0<\mu<c<1$, separately, see~\S\ref{SectionEmbedddings}.
\end{remark}

\subsection{The stratification of $\tilde{\jj}_{\mu,c}$}\label{SectionStratification}

Adapting the gluing techniques used in McDuff~\cite{MD-Stratification}, one can show that the decomposition~(\ref{LaStratification}) is a genuine stratification. In order to state its main properties, we now briefly recall the projective realizations of the ruled manifolds $\Mui$ as Hirzebruch surfaces. 

For any $\mu > 0$ and any integer $i \geq 0$ satisfying $\mu - \frac{i}{2} > 0$, let $\CP^1 \times\CP^2$ be endowed with the K\"ahler form $(\mu-a(i))\sigma_1 + \sigma_2$ where $\sigma_{j}$ is the Fubini-Study form on $\CP^{j}$ normalized so that the area of the linear $\CP^1$'s is equal to $1$, and where $a(i)$ is $\frac{i}{2}$ when $i$ is even or $\frac{i-1}{2}$ when $i$ is odd. Let $W_i$ be the Hirzebruch surface defined by setting
$$
W_i=\{([z_0,z_1],[w_0,w_1,w_2])\in\CP^1\times\CP^2 ~|~ z_0^iw_1=z_1^iw_0 \}
$$
The restriction of the projection $\pi_1:\CP^1\times\CP^2\to\CP^1$ to $W_i$ endows $W_i$ with the structure of a $\CP^1$-bundle over $\CP^1$ which is topologically $S^2 \times S^2$ if $i$ is even, and $\PbP$ if $i$ is odd. In this correspondence, the fibers are preserved and the zero section of this bundle
$$
s_0 = \{([z_0,z_1],[0,0,1])\}
$$
corresponds to a section of self-intersection $-i$ which represents the class $B^{0}-\frac{i}{2}F^{0}$ if $i$ is even or the class $B^{1} -(\frac{i-1}{2})F^{1}$ if $i$ is odd. Similarly, the section at infinity
$$s_\infty = \{([z_0,z_1],[z_{0}^{i},z_{1}^{i},1])\}$$
has self-intersection $+i$ and represents $B^{0}+\frac{i}{2}F^{0}$ or $B^{1}+\frac{i-1}{2}F^{1} $ depending on the parity of $i$ (compare with (\ref{CorrespondanceHomologique}). By the classification theorem of ruled symplectic $4$-manifolds, this correspondence establishes a symplectomorphism between $W_i$ and  $M_{\mu}^{0}$ for all even $i$'s and between $W_i$ and $M_{\mu}^{1}$ for all odd $i$'s. 

\begin{prop}[(see \cite{LP-Duke} \S4, and~\cite{MD-Stratification})]\label{StructureStratification} Let $1\leq m\leq N$. Then
\begin{enumerate}
\item The subspace $\tilde{\jj}_{\mu,c,m}$ is a smooth, co-oriented, codimension $2m$  submanifold whose closure is the union $\bigsqcup_{m\leq i\leq N}\tilde{\jj}_{\mu,c,i}$.
\item $\tilde{\jj}_{\mu,c,m}$ contains a complex structure $\tilde{J}_{m}$ coming from the blow-up of the Hirzebruch surface $W_{m}$ at a point $p$ belonging to the zero section $s_{0}$.
\item The group $\Symp_{h}(\tMuco)$ of symplectomorphisms acting trivially on homology acts smoothly on $\tilde{\jj}_{\mu,c,m}$. The stabilizer of $\tilde{J}_{m}$ is the $2$-torus $\tilde{T}_{m}$ which is generated by the lifts of the K\"ahlerian isometries of $W_{m}$ fixing $p$.
\end{enumerate}
Moreover, the decomposition $\tilde{\jj}_{\mu,c,0}\sqcup\cdots\sqcup\tilde{\jj}_{\mu,c,N}$ is a genuine stratification. Each $\tilde{\jj}_{\mu,c,m}$ has a neighborhood $\Nn_{m} \subset \tilde{\jj}_{\mu,c}$ which, once given the induced stratification, has the structure of a locally trivial fiber bundle whose typical fiber is a cone over a finite dimensional stratified space called the {\em link} of $\tilde{\jj}_{\mu,c.m}$ in $\tilde{\jj}_{\mu,c}$. Finally, the link $\Ll_{m+1}$ of $\tilde{\jj}_{\mu,c,m+1}$ in $\tilde{\jj}_{\mu,c,m}$ is a circle.
\end{prop}

Using the classifications of Hamiltonian $S^{1}$- and $T^{2}$-actions on symplectic $4$-manifolds given by Karshon~\cite{Ka-Classification} and Delzant~\cite{De}, one can show that the set of strata is in bijection with the set of conjugacy classes of maximal tori $\tilde{T}_{i}$ in $\Symp(\tilde{M}^{i}_{\mu,c})$. Finally, note that the number of strata, and hence the number of inequivalent $T^{2}$-actions, changes exactly when $\mu$ crosses an integer or when the capacity $c$ reaches the critical value $\lambda$.

\section{Homotopy type of symplectomorphism groups}\label{SectionStabilite}

\subsection{Connexity of $\Symp(\tMuci)$ and $\Symp(\tMuci,\Sigma)$}\label{SubSectionConnexite}

Let $\Symp_{h}(\tMuci)$ be the group of symplectomorphisms of $\tMuci$ acting trivially in homology. We begin this section by recalling the main result of Section 4.1 in~\cite{LP-Duke}.
\begin{prop}[(see~\cite{LP-Duke}, \S4.1)]\label{PropHomogeneousStrata}
In the stratification~(\ref{LaStratification}), the stratum $\tilde{\jj}_{\mu,c,m}$ is homotopy equivalent to the homogeneous space $\Symp_{h}(\tMuco)/\tilde{T}_{m}$.
\end{prop}

The proof of this proposition is done in three steps. First, one shows that the stratum $\tilde{\jj}_{\mu,c,m}$ is homotopy equivalent to the space $\Cc_{m}$ of all configurations of three symplectic embedded surfaces in $\tMuco$ which intersect tranversally and positively, and lie in classes $E$, $F-E$, and $D_{-m}$. Second, one proves that the configuration space $\Cc_{m}$ retracts onto the subspace $\Cc^{0}_{m}$ of configurations made of surfaces intersecting orthogonally. Finally, one shows that the symplectomorphism group $\Symp_{h}(\tMuco)$ acts transitively on $\Cc_{m}^{0}$ and that the stabilizer of a standard configuration $\Gamma_{0}$ retracts onto the $2$-torus $\tilde{T}_{m}$ of K\"ahlerian isometries of the Hirzebruch complex structure $\tilde{J}_{m}$. Note that this last step is rather subtle since it depends on the existence of a globally contracting Liouville vector field on the complement $\tMuco \setminus \Gamma_{0}$, and on Gromov's result on the contractibility of the group of compactly supported symplectomorphisms of the standard ball $B^{4}\in\R^{4}$.

\begin{lemma}[(see~\cite{LP-Duke}, Prop. 4.15)]\label{LemmaConnexity}
The group $\Symp_{h}(\tMuci)$ is connected.
\end{lemma}
\begin{proof}
We can suppose $i=0$. The open stratum $\tilde{\jj}_{\mu,c,0}$ in the stratification~(\ref{LaStratification}) is path-connected and, by Proposition~\ref{StructureStratification}, is homotopy equivalent to the homogeneous space $\Symp_{h}(\tMuco)/\tilde{T}_{0}$.
\end{proof}
\begin{cor}
The full group $\Symp(\tMuci)$ is connected except for $\widetilde{M}_{1,c}^{0}\simeq \widetilde{M}_{1-c,1-c}^{1}$ in which case it has exactly two components.
\begin{proof}
The last term in the short exact sequence $$1\to\Symp_{h}(\tMuci)\to\Symp(\tMuci)\to\Aut_{c_{1},\om}(H_{2}(\tMuci;\Z))\to 1 $$ is isomorphic to $\Z_{2}$ when $\om(F-E)=\om(B-E)$ and is trivial otherwise.
\end{proof}
\end{cor}
\begin{proof}[Proof of Lemma~\ref{RetractionSymp}]
Suppose $c\leq\mu$. According to Lemma~\ref{PersistanceEetF-E}, given any tame $J\in\jj(\tMuci)$, the exceptional class $[S]$ is always represented by a unique embedded $J$-holomorphic sphere. It follows that the space $\Cc[S]$ of all embedded symplectic spheres in class $[S]$ is homotopy equivalent to $\jj(\tMuci)$ so that the fibration
\[
\Symp_{h}(\tMuci,S)\to\Symp_{h}(\tMuci)\to\Cc[S]
\]
has a contractible base.
\end{proof}

\subsection{Inflation and the stability theorem}\label{SubSectionStabilite}

The goal of this section is to prove Theorem~\ref{StabilityTheorem} which states that when $\mu_{j}\geq c_{j}$, $j\in\{1,2\}$, the symplectomorphism groups $\Symp(\tilde{M}^i_{\mu_{1},c_{1}})$ and $\Symp(\tilde{M}^i_{\mu_{2},c_{2}})$ have the same homotopy type whenever $\mu_{2}\in (\ell_{1},\ell_{1}+1]$, and either $c_{j}<\lambda_{j}$ or $\lambda_{j} \leq c_{j}$. In other words, when we vary the cohomology class of $\tilde{\om}_{\mu,c}^{i}$, the homotopy type of $\Symp(\tilde{M}^i_{\mu,c})$ jumps exactly when a change occurs in the geometry of the corresponding stratification $\tilde{\jj}_{\mu,c}=\tilde{\jj}_{\mu,c,0}\sqcup\cdots\sqcup\tilde{\jj}_{\mu,c,N}$. Again, we will work only with $\Muo$ and we will omit the superscript ``$0$'' on associated objects.

Let $G_{\mu}=\Symp_{h}(\Muo)$ and denote by $\tGuc$ the group of symplectomorphisms of the blow-up $\tMuco$ acting trivially on homology. Consider the natural action of the identity component of the diffeomorphism group of $\tMuco$ on the space $\Ssuc$ of symplectic forms isotopic to the blown-up form $\omuc$. By Moser's theorem, the evaluation map $\ev:\phi\mapsto\phi^*(\omuc)$ defines a fibration
\begin{equation}\label{MozerFibration}
\tGuc=\Symp_{h}(\tMuco)\cap\Diff_{0}\to\Diff_0\stackrel{\ev}{\to}\Ssuc
\end{equation}
Therefore, to prove that the symplectomorphism groups $G^i_{\mu_{1},c_{1}}$ and $G^i_{\mu_{2},c_{2}}$ are homotopy equivalent, it would be enough to construct a homotopy equivalence between the spaces $\Ss_{\mu_{1},c_{1}}$ and $\Ss_{\mu_{2},c_{2}}$ such that the diagram
\begin{eqnarray}\label{FibrationSympDiffForms}
\begin{array}{ccccc}
G_{\mu_{1},c_{1}} & \to & \Diff_0   & \to & \Ss_{\mu_{1},c_{1}} \\
                    &     & \parallel &     & \updownarrow          \\
G_{\mu_{2},c_{2}} & \to & \Diff_0   & \to & \Ss_{\mu_{2},c_{2}} 
\end{array}
\end{eqnarray}
commutes up to homotopy. Since we do not know any way of constructing such a map directly, we will proceed a little differently. Indeed, we will use a nice idea of McDuff that consists in considering the action of $\Diff_0$ on the larger space $\Xx_{\mu_{i},c_{i}}$ of pairs $(\om, J)$ such that $\om\in\Ss_{\mu_{i},c_{i}}$ and $J$ is tamed by $\om$. This space projects to $\Ss_{\mu_{i},c_{i}}$ and to the space $\Aa_{\mu_{i},c_{i}}$ of almost complex structures tamed by some element of $\Ss_{\mu_{i},c_{i}}$. It turns out that these projections are Serre fibrations with contractible fibers. This implies that the three spaces $\Ss_{\mu_{i},c_{i}}$, $\Xx_{\mu_{i},c_{i}}$ and $\Aa_{\mu_{i},c_{i}}$ are homotopy equivalent. Moreover, the previous fibrations factorize through the spaces $\Xx_{\mu_{1},c_{1}}$ and $\Xx_{\mu_{2},c_{2}}$ as in the following diagram:
\begin{eqnarray}\label{ExtendedFibration}
\begin{array}{ccccccc}
             &     &           &   &\Aa_{\mu_{1},c_{1}}&   & \\

             &     &           &   &\uparrow             &   & \\

G_{\mu_{1},c_{1}}&\to&\Diff_0&\to&\Xx_{\mu_{1},c_{1}}&\to&\Ss_{\mu_{1},c_{1}}\\

             &     & \parallel &   &                     &   & \\

G_{\mu_{2},c_{2}}&\to&\Diff_0&\to&\Xx_{\mu_{2},c_{2}}&\to&\Ss_{\mu_{2},c_{2}}\\

             &     &           &   &\downarrow           &   & \\

             &     &           &   &\Aa_{\mu_{2},c_{2}}&   & \\
\end{array}
\end{eqnarray}

Therefore, to show that the groups $G_{\mu_{1},c_{1}}$ and $G_{\mu_{2},c_{2}}$ are homotopy equivalent, it is sufficient to find a homotopy equivalence $\Aa_{\mu_{1},c_{1}}\to \Aa_{\mu_{2},c_{2}}$ that commutes, up to homotopy, with the action of $\Diff_0$. Indeed, we prove a stronger statement, namely:

\begin{prop}\label{Egalite}
Let $\mu_{1}=\ell+\lambda_{1}\geq 1$. Given $\mu_{2}\in (\ell,\ell+1]$, write $\mu_{2}=\ell+\lambda_{2}$. Consider $c_{1},c_{2}\in(0,1)$ such that either $c_{1}<\lambda_{1}$ and $c_{2}<\lambda_{2}$, or $\lambda_{1}\leq c_{1}$ and $\lambda_{2}\leq c_{2}$. Then the spaces $\Aa_{\mu_{1},c_{1}}$ and $ \Aa_{\mu_{2},c_{2}}$ are the same.
\end{prop} 

The proof is based on the Inflation Lemma of Lalonde-McDuff which asserts that, in the presence of a $J$-holomorphic submanifold with nonnegative self-intersection,  it is possible to deform a symplectic form $\om_{0}$ through a family of forms $\om_{t}$ tamed by the almost complex structure $J$. 
\begin{prop*}[(Inflation lemma, see~\cite{La,MD-Haefliger})]\label{Inflation}
Let $J\in\J(\om_{0})$ be a tame almost complex structure on a symplectic $4$-manifold $(M,\om_{0})$ that admits a $J$-holomorphic curve $Z$ with $Z\cdot Z\geq 0$. Then there is a family
$\om_{t}$, $t\geq 0$, of symplectic forms that all tame $J$ and have cohomology class $[\om_{t}]=[\om_{0}]+t\PD(Z)$,
where $\PD(Z)$ denotes the Poincar\'e dual of the homology class $[Z]$.
\end{prop*}
\begin{proof}
See Lemma 3.1 in~\cite{MD-Haefliger}.
\end{proof}

\begin{proof}[Proof of Proposition~\ref{Egalite}] 
The proof is done in two steps. We first prove that for $\mu=\ell+\lambda\geq 1$ and $c<\lambda$ (resp. $\lambda\leq c$), the spaces $\Aa_{\mu,c}$ and $\Aa_{\mu',c}$ are the same, for all $\mu'\in(\ell+c,\mu]$ (resp. $\mu'\in(\ell,\mu])$. In the second step, we show that $\Aa_{\mu,c}=\Aa_{\mu,c'}$ whenever the two parameters $c$ and $c'$ verify $c\leq c'<\lambda$ or $\lambda\leq c\leq c'$.

{\it Step 1}: We begin by showing that for any $\mu\geq 1$ and any $\epsilon>0$, the space $\Aa_{\mu,c}$ is included in $\Aa_{\mu+\epsilon,c}$. Indeed, Lemma~\ref{ModuliSpaceFiber} implies that for all $(\om_{\mu,c},J)\in\Xx_{\mu,c}$, the fiber class $F$ is represented by some embedded $J$-holomorphic sphere $S$. Therefore, it is always possible to inflate the form $\om_{\mu,c}$ along $S$ to get a one-parameter family of symplectic forms in classes $[\om_{\mu,c}]+\epsilon\PD(F)=\om_{\mu+\epsilon,c}$. Since all these forms tame $J$, this prove that $\Aa_{\mu,c}\subset \Aa_{\mu+\epsilon,c}$.

We next show that for $\mu=\ell+\lambda\geq 1$ and $c<\lambda$, the space $\Aa_{\mu,c}$ in included in $\Aa_{\mu',c}$ whenever $\mu'\in(\ell+c,\mu]$. This is done by inflating $\om_{\mu,c}$ along curves representing the classes $D_{2\ell+1}=B+\ell F$ and $D_{2\ell+2}=B+(\ell+1)F-E$. Note that when $c<\lambda$, Lemma~\ref{ExistenceCourbesDi} implies that these two classes are always represented by embedded $J$-holomorphic spheres. Writing $d_{i}$ for an appropriate $2$-form representing $\PD(D_{i})$, we define $\om_{ab}$ for $a,b\geq 0$ by 
$$\om_{ab}=\frac{\left(\om_{\mu,c} + ad_{2\ell+1}+bd_{2\ell+2}\right)}{1+a+b}$$
We then have $\om_{ab}(F)=1$, $\om_{ab}(B)=(\mu+a\ell+b(\ell+1))/(1+a+b)$, and $\om_{ab}(E)=(c+b)/(1+a+b)$. The area of the exceptional class $E$ can be made equal to $c$ by setting $b=ac/(1-c)\geq 0$. This gives us a one-parameter family of forms $\om_{a}$ verifying $\om_{a}(F)=1$, $\om_{a}(B)=((1-c)\mu+(\ell+c)a)/(1-c+a)$, and $\om_{a}(E)=c$. Letting $a\to\infty$, this shows that the almost complex structure $J$ is tamed by a symplectic form in class $[\om_{\ell+c+\epsilon,c}]$, for all $0<\epsilon\leq \lambda-c$. Therefore $\Aa_{\mu,c}\subset\Aa_{\mu',c}$, for all $\mu'\in(\ell+c,\mu]$.

Now, let prove that for $\mu=\ell+\lambda\geq 1$ and $\lambda\leq c$, the space $\Aa_{\mu,c}$ in included in $\Aa_{\mu',c}$ whenever $\mu'\in(\ell,\mu]$. Since we assume $\lambda\leq c$, the classes $D_{2\ell}$ and $D_{2\ell+1}$ are both represented by embedded $J$-holomorphic spheres, for any choice of $\om$-tamed $J$. By inflating $\om_{\mu,c}$ along these curves, we get a two-parameter family of symplectic forms $\om_{ab}=(\om_{\mu,c} + ad_{2\ell+1}+bd_{2\ell})/(1+a+b)$ for which $\om_{ab}(F)=1$, $\om_{ab}(B)=(\mu+a\ell+b\ell)/(1+a+b)$, and $\om_{ab}(E)=(c+b)/(1+a+b)$. Again, we make the area of the exceptional class $E$ equal to $c$ by setting $b=ac/(1-c)\geq 0$. We then get a one-parameter family of forms $\om_{a}$ verifying $\om_{a}(F)=1$, $\om_{a}(B)=((1-c)\mu+a\ell)/(1-c+a)$, and $\om_{a}(E)=c$. Letting $a\to\infty$, this shows that $\Aa_{\mu,c}\subset\Aa_{\mu',c}$ for all $\mu'\in(\ell,\mu]$.

{\it Step 2}: 
Let us show that $\Aa_{\mu,c'}\subset\Aa_{\mu,c}$ whenever $c<c'$. Given $J\in\Aa_{\mu,c'}$, this amounts to constructing a symplectic form that tames $J$ and for which $B$, $F$ and $E$ have area $\mu$, $1$ and $c$ respectively. By lemma~\ref{ExistenceCourbesDi} and lemma~\ref{ModuliSpaceFiber}, there exist embedded $J$-holomorphic spheres representing the classes $D_{2\ell+1}=B+\ell F$ and $F$. The inflation lemma implies that the Poincar\'e duals of $D_{2\ell+1}$ and $F$ are represented by $2$-forms $d_{2\ell+1}$ and $f$ such that the form $\om_{ab} = (\om_{\mu,c'} + ad_{2\ell+1} + bf)/(1+a)$ is symplectic for all $a,b\geq 0$. We then have $\om_{ab}(F)=1$ and $\om_{ab}(B)=(\mu+a\ell + b)/(1+a)$. Setting $b=(\mu-\ell)a=\lambda a\geq 0$, the resulting family of symplectic forms $\om_{a}$ verify $\om_{a}(F)=1$, $\om_{a}(B)=\mu$, and $\om_a(E)=c'/(1+a)$. Since $a$ can take any non negative value, this implies that $J\in\Aa_{\mu,c}$, for all $c<c'$.

We now suppose that $0<c<\lambda$ and prove that given any real number $c'\in(c,\lambda)$, the space $\Aauc$ is included in $\Aa_{\mu,c'}$. This is done by inflating along embedded $J$-holomorphic spheres in classes $D_{2\ell+2}=B+(\ell+1)F-E$ and $D_{2\ell+1}=B+\ell F$. Indeed, setting $\om_{ab} = (\om_{\mu,c} + ad_{2\ell+2} + bd_{2\ell+1})/(1+a+b)$ where $b=(\ell+1-\mu)/(\mu-\ell) \geq 0$, the resulting one-parameter family of symplectic forms verify $\om_a(F)=1$, $\om_a(B)=\mu$, and $\om_a(E)=\lambda(a+c)/(\lambda + a)$. This shows that any almost complex structure $J\in\Aa_{\mu,c}$ also belongs to $\Aa_{\mu,c'}$

Finally, given $\lambda \leq c < c'$, we need to prove that $\Aa_{\mu,c}\subset\Aa_{\mu,c'}$. Since $\lambda \leq c$, Lemma~\ref{ExistenceCourbesDi} implies that the classes $D_{2\ell}=B+\ell F-E$ and $F$ are both represented by embedded $J$-holomorphic spheres, for all $J$. We define $\om_{ab}=(\om_{\mu,c}+ad_{2\ell}+bf)/(1+a)$ and set $b=a(\mu-\ell)\geq 0$ to get symplectic forms $\om_{a}$, $a\geq 0$, such that $\om_a(B)=\mu$, $\om_a(F)=1$ and $\om_a(E)=(a+c)/(a+1)$. Therefore, $J\in\Aa_{\mu,c'}$, for all $\lambda\leq c<c'<1$. This completes the proof of Proposition~\ref{Egalite}.
\end{proof}

\begin{remark}
Proposition~\ref{Egalite} is optimal in the sense that it cannot be improved by performing the inflation process along curves representing other classes.
\end{remark}

We can now prove the stability theorem.

\begin{proof}[Proof of Theorem~\ref{StabilityTheorem}]
Consider the diagram~\ref{ExtendedFibration}. By Proposition~\ref{Egalite}, we know that the spaces $\Aauc$ and $\Aa_{\mu,c'}$ are equal whenever  $c < c' < \lambda$ or $\lambda \leq c < c'$. Therefore, the spaces of symplectic forms $\Ssuc$ and $\Ss_{\mu,c'}$ are homotopy equivalent and this implies that the group $\tGuc$ has the same homotopy type as $\tilde{G}_{\mu,c'}$. Finally, in the special case $\mu=1$, the full symplectomorphism group of the blow-up is a $\Z_{2}$-extension of its identity component and is homotopy equivalent to $T^{2}\times\Z_{2}$ for all $0<c<1$.
\end{proof}

We note that the stratification $\tilde{\jj}_{\mu,c}=\tilde{\jj}_{\mu,c,0}\sqcup\cdots\sqcup\tilde{\jj}_{\mu,c,N}$ induces similar decompositions of the spaces $\Xx_{\mu,c}$ and $\Aa_{\mu,c}$ into Fr\'echet submanifolds invariant under the action of $\Diff_{0}$. The proof of Proposition~\ref{Egalite} shows that, for all $\mu\geq 1$, $c\in(0,1)$, and $\epsilon>0$, we have an inclusion $\Aa_{\mu,c}\into \Aa_{\mu+\epsilon,c}$. Moreover, if we increase $\mu$ while keeping $c$ constant, we see that the space $\Aa_{\mu,c}$ changes by the addition of a stratum of codimension $4n-2$ each time $\mu$ crosses an integer $n$, or by the addition of a stratum of codimension $4n$ each time $\mu$ crosses a real number of the form $n+c$. Consequently, the inclusion $\Aa_{\mu,c}\into \Aa_{\mu+\epsilon,c}$ is $(4\ell-2)$-connected if $c\geq\lambda$, and $4\ell$-connected if $c<\lambda$. At the level of symplectomorphism groups, the homotopy commuting diagram~(\ref{ExtendedFibration}) implies that the the family of inclusions $\Aa_{\mu,c}\into \Aa_{\mu+\epsilon,c}$ induces homotopy coherent diagrams 
\[
\begin{array}{rcl}
G_{\mu,c} & \to & G_{\mu+\epsilon,c}\\
& \searrow & \downarrow \\
& & G_{\mu+\epsilon+\delta,c}
\end{array}
\]
in which the map $G_{\mu,c}\into G_{\mu+\epsilon,c}$ is $(4\ell-3)$-connected if $c\geq\lambda$, and $(4\ell-1)$-connected if $c<\lambda$. In particular,

\begin{prop}\label{IsomorphismFundamentalGroups}
For each $c\in(0,1)$ and $\epsilon\in(0,c)$, the map $G_{1+\epsilon,c}\to G_{\mu,c}$ induces an isomorphism of fundamental groups.
\end{prop}\qed

\subsection{Homotopy equivalences with stabilizer subgroups}\label{SubSectionHomotopySymp}

We now prove Theorem~\ref{HomotopyOfSymp} which states that the symplectomorphism group $\Symp(\tMuco)$ is homotopy equivalent to the stabilizer of a point in an associated ruled manifold.

\begin{proof}[Proof of Theorem~\ref{HomotopyOfSymp}]
As before, it suffices to consider the manifold $\tMuco$ obtained by blowing up $\Muo$ at a ball $B_{c}$ of capacity $c$. First, suppose that $c<\lambda$. We want to show that the group $\Symp(\tMuco)$ is homotopy equivalent to the stabilizer of a point in $\Muo$. By Lemma~\ref{RetractionSymp} above and Lemma~2.3 of~\cite{LP-Duke}, we know that $\Symp(\tMuco)$ is homotopy equivalent to the group $\Symp^{\U(2)}(\Muo,B_{c})$ of those symplectomorphisms which act $\U(2)$-linearly in a neighborhood of $B_{c}$. Note that this group is included in the stabilizer of the center $p\in B_{c}$. Moreover, Theorem~\ref{StabilityTheorem} implies that for any $c'\in (0,c)$, the inclusion $i_{c,c'}:\Symp^{\U(2)}(\Muo,B_{c})\into\Symp^{\U(2)}(\Muo,B_{c'})$ is a homotopy equivalence. 
Furthermore, by Moser argument, every compact subset of the stabilizer $\Symp(\Muo,p)$ is homotopic to a family of symplectomorphisms that act linearly on a ball of sufficiently small capacity centered at $p$. Therefore, the natural map 
$$\lim_{\longrightarrow} \Symp^{\U(2)}(\Muo,B_{c}) \to \Symp(\Muo,p)$$
between the direct limit defined by the inclusions $\{i_{c,c'}\}$ and the stabilizer is a (weak) homotopy equivalence.
 
When $\lambda\leq c$, the blow-up of $\Muo$ at $B_{c}$ is symplectomorphic to the blow-up of a small ball $B_{1-c}$ in the $c$-related manifold $M_{\mu}^{0,c}=M^{1}_{\mu-c}$. The same argument as before shows that in this case $\Symp(\tMuco)\simeq \Symp(M^{1}_{\mu-c},p)$. 
\end{proof}

\subsection{The group $\Symp(\tMucl,\Sigma)$ when $0<\mu<c<1$}\label{SectionHomotopySpecialCase}

As explained in the Introduction, when $0<\mu<c<1$, the manifold $\tMucl$ is conformally symplectomorphic to $\widetilde{M}^{0}_{\mu',c'}$ where $\mu'=1/\mu+1-c$ and $c'=(1-c)/1+\mu-c$. This identification takes the class $[\Sigma]$ to $B^{0}-E^{0}$ so that the space $\jj_{\mu',c'}^{0}$ is stratified according to the degeneracy types of curves representing $[\Sigma]$. The space $\Ss([\Sigma])$ of exceptional symplectic spheres in class $[\Sigma]$ is then homotopy equivalent to the open stratum $\jj_{\mu',c',0}^{0}$ which, in turn, is homotopy equivalent to $\Symp(\tM_{\mu',c'})/T_{0}$. Looking at the fibration
\[
\Symp(\tMucl,\Sigma)\simeq\Symp(\tM_{\mu',c'}^{0},\Sigma) \to \Symp(\tM_{\mu',c'}^{0}) \to \Ss([\Sigma])\simeq \Symp(\tM_{\mu',c'}^{0})/T_{0}
\]
we conclude that $\Symp(\tMucl,\Sigma)\simeq T_{0}$.

In order to compare $\Symp(\tMucl,\Sigma)$ with the stabiliser of a point $p\in M_{\mu}^{1}$, we briefly describe a relative version of diagram~(\ref{ExtendedFibration}). In fact, because we will use the same technique in our analysis of two disjoint balls in $\CP^{2}$, we consider more generally the group $\Symp(X,\Dd)$ of symplectomorphisms which leave invariant each curve in the union $\Dd$ of $n$ disjoint, embedded, symplectic spheres representing some classes $D_{1},\ldots, D_{n}$ in a symplectic $4$-manifold $(X,\om)$. Denote by $\Ss_{0}(\om,\Dd)$ the connected component of $\om$ in the set of all symplectic forms cohomologous to $\om$ and for which $\Dd$ is symplectic. Standard arguments involving the Mozer Stability Theorem in dimension $2$, the extension of Hamiltonian isotopies, and the relative version of Moser theorem, imply that the identity component of the group $\Diff(X,\Dd)$ (made of those diffeomorphisms sending each curve of $\Dd$ to itself) acts transitively on $\Ss_{0}(\om,\Dd)$. This action defines a fibration
\[
\Symp(X,\Dd)\cap\Diff_{0}(X,\Dd)\to\Diff_{0}(X,\Dd)\to\Ss_{0}(\om,\Dd)
\]
whose base is homotopy equivalent to the space $\Xx(\om,\Dd)$ of pairs $(\tau,J)$ such that $\tau\in\Ss_{0}(\om,\Dd)$, $J$ tames $\om$, and $\Dd$ is $J$-holomorphic. As before, the projection on the second factor $\Xx(\om,\Dd)\to\Aa(\om,\Dd)$ is an homotopy equivalence so that we can compare the symplectic stabilizer of $\Dd$ associated to different symplectic forms through the diagram 
\begin{eqnarray}\label{RelativeExtendedFibration}
\begin{array}{ccccccc}
             &     &           &   &\Aa(\om,\Dd)&   & \\

             &     &           &   &\uparrow             &   & \\

\Symp(X,\Dd;\om)\cap\Diff_{0}(X,\Dd)&\to&\Diff_0&\to&\Xx(\om,\Dd)&\to&\Ss_{0}(\om,\Dd)\\

             &     & \parallel &   &                     &   & \\

\Symp(X,\Dd;\om')\cap\Diff_{0}(X,\Dd)&\to&\Diff_0&\to&\Xx(\om',\Dd)&\to&\Ss_{0}(\om',\Dd)\\

             &     &           &   &\downarrow           &   & \\

             &     &           &   &\Aa(\om',\Dd)&   & \\
\end{array}
\end{eqnarray}

\begin{proof}[Proof of Theorem~\ref{HomotopyOfSympSpecialCase}]
In the case of the (connected) group $\Symp(\tMucl,\Sigma)$ with $0<\mu<c<1$, the space $\Aa(\om,\Sigma)$ belongs to the open stratum $\Aa_{\mu',c',0}\subset\Aa_{\mu',c'}$. Therefore, for each $J\in\Aa(\om,\Sigma)$, there exist embedded $J$-holomorphic spheres in classes $F^{0}$, $B^{0}$, and $B^{0}+F^{0}-E^{0}$. Performing the inflation along theses curves (as in the proof of Theorem~\ref{StabilityTheorem}) shows that $\Aa(\om,\Sigma)$ in independent of $\mu$ and $c$,   so that we get a coherent family of homotopy equivalences $\Symp(\tMucl,\Sigma)\to\Symp(\tM_{\mu'',c''},\Sigma)$ whenever $0<\mu''<c''<1$. Because the group $\Symp(M^{1}_{\mu})$ is homotopy equivalent to $\U(2)$ when $\mu\leq 1$ (see Theorem~\ref{Gromov} below), we can look at the direct limit of the groups $\Symp(\tM_{\mu=c/2,c},\Sigma)$ as $c\to 0$ to conclude that $\Symp(\tMucl,\Sigma)$ is homotopy equivalent to the stabilizer of a point in $\Mul$ when $0<\mu<c<1$.
\end{proof}

\subsection{Formal consequences}\label{SectionFormalConsequences}

The homotopy type of the groups $\Symp(\Mui)$ was first investigated by Gromov in his celebrated paper~\cite{Gr} in which he proved the following results using pseudo-holomorphic techniques and cut-and-paste arguments:

\begin{thm}[(see Gromov~\cite{Gr})]\label{Gromov}
\begin{enumerate}
\item When $\mu=1$, the symplectomorphism group of the manifold $\Muo$ retracts onto its subgroup $(\SO(3)\times\SO(3))\ltimes\Z_2$ of k\"ahlerian isometries.
\item  When $0<\mu\leq 1$, the symplectomorphism group of the manifold $\Mul$ retracts onto its subgroup $\U(2)$ of k\"ahlerian isometries.
\end{enumerate}
\end{thm}

Later, Abreu~\cite{Ab} and then Abreu-McDuff~\cite{AM} extended Gromov's techniques and computed the rational cohomology ring $H^{*}(\Symp(\Mui);\Q)$. The rational cohomology ring $H^{*}(\BSymp(\Mui);\Q)$ was computed by Abreu-Graja-Kitchloo~\cite{AGK}. Their results can be summarized as follow:

\begin{thm}[(see~\cite{AM})]\label{AbreuMcDuff}
When $\mu>1$, the rational cohomology ring of the group $\Symp(\Mui)$ is isomorphic to the product $\Lambda(t,x,y)\otimes \Q[w_\ell]$ with generators $t$, $x$, $y$ and $w_\ell$ of degree $1$, $3$, $3$ and~$4\ell+2i$. Therefore, the rational homotopy groups of $\Symp(\Mui)$ vanish in all dimensions except in dimensions $1$, $3$, and~$4\ell+2i$, in which cases we have $\pi_{1}=\Q$, $\pi_{3}=\Q^{2}$, and $\pi_{4\ell+2i}=\Q$.  
\end{thm}

\begin{thm}[(see~\cite{AM} and~\cite{AGK})]\label{AbreuMcDuff-Classifying}
When $\mu>1$, the rational cohomology ring of the classifying space $\BSymp(\Mui)$ is isomorphic to
$$\Q[T,X,Y]/ R_{\ell}$$
where $R_{\ell}$ is the ideal generated by 
\begin{align}
& T(T^{2}+X-Y)(T^{2}+16X-4Y)\cdots (T^{2}+\ell^{4}X-\ell^{2}Y)&\text{in the case $i=0$,}\\
& (T^{2}+X-Y)(T^{2}+3^{4}X-3^{2}Y)\cdots (T^{2}+(2\ell+1)^{4}X-(2\ell+1)^{2}Y)&\text{in the case $i=1$,}
\end{align}
and where the generators $T$, $X$, $Y$ are of degree $2$, $4$, and~$4$.
\end{thm}

Combining these results with Theorem~\ref{HomotopyOfSymp}, we immediately get the following corollaries:

\begin{cor} \label{RationalHomotopy} 
Let $\tMuci$ be the symplectic blow-up of $\Mui$ at a ball of capacity $c\in(0,1)$ and suppose that $\mu>1$. Then,
\begin{itemize}
\item[1)] when $c<\lambda$, the non trivial rational homotopy groups of the topological group $\Symp(\tMuci)$ are $\pi_{1}\simeq\Q^{3}$ and \mbox{$\pi_{4\ell+2i}\simeq\Q$}.

\item[2)] when $\lambda\leq c$, the non trivial rational homotopy groups of the topological group $\Symp(\tMuci)$ are $\pi_{1}\simeq\Q^{3}$ and $\pi_{4\ell+2i-2}\simeq\Q$.
\end{itemize}
\end{cor}

\begin{proof}
The rational homotopy groups of $\Symp(\tMuci)$ are easily computed by combining Theorem~\ref{HomotopyOfSymp} with Theorem~\ref{AbreuMcDuff}. Indeed, suppose we blow-up a small ball in $\Muo$, where $\mu=\ell+\lambda>1$. Then we have a homotopy fibration
$$\Symp(\tMuco)\to\Symp(\Muo)\to S^{2}\times S^{2}$$
which gives an exact homotopy sequence whose only non trivial terms are
$$0\to\pi_{4\ell}(\Symp(\tMuco))\to\Q\to 0$$
and
\begin{multline}
0\to\pi_{3}(\Symp(\tMuco))\to\Q^{2}\to\Q^{2}\to\pi_{2}(\Symp(\tMuco))\to\\
0\to\Q^{2}\to\pi_{1}(\Symp(\tMuco))\to\Q\to 0
\end{multline}
In these exact sequences, the only problem is to understand the map $\partial:\Q^{2}\simeq\pi_{3}(\Symp(\Muo))\to\pi_{2}(S^{2}\times S^{2})\simeq \Q^{2}$. But from the work of Abreu-McDuff, the group $\pi_{3}(\Muo)\otimes\Q$ is identified with $\pi_{3}(\SO(3)\times\SO(3))\otimes\Q$ under the isometric action of $\SO(3)\times\SO(3)$ on $S^{2}\times S^{2}$. Therefore, the map $\partial$ is an isomorphism and both $\pi_{3}(\Symp(\tMuco))$ and $\pi_{2}(\Symp(\tMuco))$ are trivial.

When we blow-up a ball of capacity $c\geq\lambda$ in $\Muo$, the group $\Symp(\tMuco)$ is homotopy equivalent to the stabilizer of a point in $\tilde{M}^{1}_{\mu-c}$. Since $\lfloor\mu-c\rfloor = \lfloor\ell+\lambda-c\rfloor=\ell-1$, the only non trivial homotopy groups of $\Symp(\tilde{M}^{1}_{\mu-c})$ are in degree $1$, $3$ and $4(\ell-1)+2=4\ell-2$. When $\ell\geq 2$, the homotopy fibration
$$\Symp(\tMuco)\to\Symp(\tilde{M}^{1}_{\mu-c})\to \NTB$$
gives rise to an exact homotopy sequence whose only non trivial terms are
$$0\to\pi_{4\ell-2}(\Symp(\tMuco))\to\Q\to 0$$
and
\begin{multline}
0\to\pi_{3}(\Symp(\tMuco))\to\Q^{2}\to\Q^{2}\to\pi_{2}(\Symp(\tMuco))\to\\ 
0\to\Q^{2}\to\pi_{1}(\Symp(\tMuco))\to\Q\to 0
\end{multline}
When $\ell=1$, the exact homotopy sequence reduces to
\begin{multline}
0\to\pi_{3}(\Symp(\tMuco))\to\Q^{2}\to\Q^{2}\to\pi_{2}(\Symp(\tMuco))\to\\ \pi_{2}(\Symp(M_{\mu-c}^{1})
\stackrel{0}{\to}\Q^{2}\to\pi_{1}(\Symp(\tMuco))\to\Q\to 0
\end{multline}
In both cases, Lemma 2.11 in Abreu-McDuff~\cite{AM} implies that the homotopy group $\pi_{3}(\Symp(\tilde{M}^{1}_{\mu-c}))$ is isomorphic to $\pi_{3}(\Dd)\otimes\Q\simeq\Q^{2}$, where $\Dd$ is the group of orientation preserving fiberwise diffeomorphisms of $\NTB$, and that the connecting map $\partial:\pi_{3}(\Symp(\tilde{M}^{1}_{\mu-c}))\to\pi_{2}(\NTB)\simeq \Q^{2}$ is an isomorphism. Therefore, $\pi_{3}(\Symp(\tMuco))$ is trivial while $\pi_{2}(\Symp(\tMuco))\simeq\pi_{2}(\Symp(M_{\mu-c}^{1})$.
\end{proof}

Note that the vector space $\pi_{*}(\Symp(\tMuci))\otimes\Q$ has $3$ \emph{topological} generators of degree $1$ which do not depend on the cohomology class of the symplectic form, and one \emph{symplectic} generator whose degree changes precisely when $\mu$ crosses an integer or when the capacity $c$ reaches the critical value $\lambda$. The topological generators are \emph{robust} in the sense that their images in $\pi_{*}(\Diff)\otimes\Q$ generate a $3$ dimensional subspace, while the symplectic generator is \emph{fragile}, that is, it does not survive the inclusion $\Symp\into\Diff$.

\begin{cor} \label{RationalCohomology} 
Let $\tMuci$ be the symplectic blow-up of $\Mui$ at a ball of capacity $c\in(0,1)$. In the case $i=0$, suppose that $\mu>1$, while in the case $i=1$ suppose that $c\leq\mu$. Then,
\begin{itemize}
\item[1)] when $c<\lambda$, the rational cohomology ring of the group $\Symp(\tMuci)$ is isomorphic to the product $\Lambda(\alpha_1,\alpha_2,\alpha_3)\otimes\Q[\epsilon]$ with generators $\alpha_1, \alpha_2, \alpha_3,$ and $\epsilon$ of degree $1,1,1,$ and $4\ell + 2i$.
\item[2)] When $\lambda\leq c$, the rational cohomology ring of the group $\Symp(\tMuci)$ is isomorphic to the product $\Lambda(\alpha_1,\alpha_2,\alpha_3)\otimes\Q[\epsilon]$ with generators $\alpha_1, \alpha_2, \alpha_3,$ and $\epsilon$ of degree $1,1,1,$ and $4\ell+2i-2$.
\end{itemize}
\end{cor}
\begin{proof}
Let $G$ be a connected topological group. The Samelson product 
$\langle\cdot,\cdot\rangle:\pi_{p}(G)\otimes\pi_{q}(G) \to \pi_{p+q}(G)$ defined by the commutator map 
$$
\begin{array}{cccl}
\langle\alpha,\beta\rangle: & S^{p}\wedge S^{q} & \to & G\\
 & (s,t) &\mapsto& \alpha(s)\beta(t)\alpha^{-1}(s)\beta^{-1}(t)
\end{array}
$$
gives the rational homotopy of $G$ the structure of a graded Lie algebra. The rational homology $H_{*}(G;\Q)$ is a connected and cocommutative graded Hopf algebra when equipped with the Pontrjagin product induced by the multiplication map. The Milnor-Moore theorem states that the Hurewicz homomorphism $\pi_{*}(G)\otimes\Q\to H_{*}(G;\Q)$ induces an isomorphism of Hopf algebras between the universal enveloping algebra $\Uu(\pi_{*}(G)\otimes\Q)$ and $H_{*}(G;\Q)$. In particular, it implies that $H_{*}(G;\Q)$ is always generated by its spherical classes.

In our situation, Corollary~\ref{RationalHomotopy} implies that the rational homology of $\Symp(\tMuci)$ is of finite type. Consequently, the rational cohomology $H^{*}(\Symp(\tMuci;\Q)$ is a commutative Hopf algebra dual to the homology Hopf algebra. By the Borel classification of commutative Hopf algebra over fields of characteristic zero, this implies that
$$H^{*}(\Symp(\tMuci;\Q)\simeq \Lambda(\alpha_1,\alpha_2,\alpha_3)\otimes\Q[\epsilon]$$ 
where $\alpha_1, \alpha_2, \alpha_3$ are three generators of degree $1$, and $\epsilon$ is a generator of degree $4\ell + 2i$ or $4\ell + 2i-2$ depending on whether $c<\lambda$ or $c\geq\lambda$.
\end{proof}

Since the group $\Symp(\Mui)$ acts transitively on $\Mui$, the classifying space of the isotropy group $\Symp(\Mui,p)$ can be understood from the homotopy fibration
$$\Mui\to B\Symp(\Mui,p)\to B\Symp(\Mui)$$ 
This leads to the following two corollaries:

\begin{cor} \label{CohomologyClassifyingSpace} 
Let $\tMuci$ be as in Corollary~\ref{RationalHomotopy}. Then,
\begin{itemize}
\item[1)] when $c<\lambda$, the rational cohomology module of the classifying space $\BSymp(\tMuci)$ is isomorphic to $H^{*}(\BSymp(\Mui))\otimes H^{*}(\Mui)$.
\item[2)]  When $\lambda\leq c$, the rational cohomology mudule of the classifying space $\BSymp(\tMuci)$ is isomorphic to $H^{*}(\BSymp(\Mdual))\otimes H^{*}(\Mdual)$.
\end{itemize}
\end{cor}

\begin{proof}
By Theorem~\ref{AbreuMcDuff-Classifying}, the rational cohomology rings of both $\Mui$ and $B\Symp(\Mui)$ are supported in even degrees. Consequently, the spectral sequence associated to this fibration collapses at the second stage $E_{2}^{p,q}$. Therefore,
$$H^{*}(B\Symp(\Mui,p);\Q) \simeq H^{*}(\Mui;\Q)\otimes H^{*}(B\Symp(\Mui);\Q)$$
and the result follows from Theorem~\ref{HomotopyOfSymp}.
\end{proof}

\begin{cor} \label{SamelsonProducts} The Pontrjagin product on $H_{*}(\Symp(\tMuco);\Q)$ is (graded) commutative and the Samelson product on $\pi_{*}(\Symp(\tMuco))\otimes\Q$ is trivial, except when $1<\mu<2$ and $\lambda\leq c$ in which case the Samelson products $\langle \alpha_{1},\alpha_{3}\rangle=-\langle \alpha_{2},\alpha_{3}\rangle$ are nonzero multiples of $\epsilon\in\pi_{2}(\Symp(\tMuco))\otimes\Q$.  
\end{cor}
\begin{proof}
(See~\cite{LP-Duke} for an alternative proof). For dimensional reasons, the Samelson products defined on $\pi_{*}(\Symp(\tMuco))\otimes\Q$ are all trivial except possibly when $1<\mu<2$ and $\lambda\leq c$. In that case, the jumping generator $\epsilon$ has degree $2$ and hence the products $\langle \alpha_{i},\alpha_{j}\rangle$ might be nonzero. By~\cite{LP-Duke} Proposition~4.7, we can choose the classes $\alpha_{1}, \alpha_{2}, \alpha_{3}$ such that $\alpha_{1},\alpha_{2}$ correspond to lifts (to the blow-up) of the generators of the standard K\"ahlerian (and Hamiltonian) $T^{2}$-action on the Hirzebruch surface $W^{1}$, while $\alpha_{1}+\alpha_{2},\alpha_{3}$ correspond to lifts of the generators of the K\"ahlerian $T^{2}$-action on 
$W^{2}$. Consequently, $\langle\alpha_{1},\alpha_{2}\rangle=0$ while $\langle\alpha_{1},\alpha_{3}\rangle = - \langle\alpha_{2},\alpha_{3}\rangle$.

Now, to understand the product $\langle\alpha_{1},\alpha_{3}\rangle$, we use the fact that it corresponds (up to sign) to a rational Whitehead product in the homotopy of the classifying space $B\Symp(\tMuco)$. Since the non trivial homotopy groups of $B\Symp(\tMuco)$ are $\pi_{2}\simeq\Q^{3}$ and $\pi_{3}\simeq\Q$, the minimal model of $B\Symp(\tMuco)$ has only four generators $a_{2},b_{2},c_{2},e_{3}$ (here the index is the degree). By the Hurewitz theorem, $\pi_{2}\simeq H_{2}\simeq H^{2}$ which implies that $a_{2},b_{2},c_{2}$ are cocycles. Moreover, by Corollary~\ref{CohomologyClassifyingSpace} we have
$$H_{k}(B\Symp(\tMuco))\simeq
\begin{cases}
0, & \text{if $k$ is odd,}\\
\Q^{k+1}, & \text{if $k$ is even}
\end{cases} $$
and this implies that $e_{3}$ cannot be a cocycle. Thus 
$$de_{3}=\text{~a linear combination of $a^{2}_{2},b^{2}_{2},c^{2}_{2}, a_{2}b_{2}, a_{2}c_{2}$, and $b_{2}c_{2}$}$$
By the description of the rational Whitehead products given by Andrews-Arkowitz~\cite{AA}, this implies that the Whitehead product map 
$$[-,-]:\pi_{2}(B\Symp(\tMuco))\otimes\pi_{2}(B\Symp(\tMuco))\to\pi_{3}(B\Symp(\tMuco))$$
is non trivial or, after desuspension, that the Samelson product map
$$\langle -,-\rangle:\pi_{1}(\Symp(\tMuco))\otimes\pi_{1}(\Symp(\tMuco))\to\pi_{2}(\Symp(\tMuco))$$ 
is non trivial, that is, the product $\langle \alpha_{1},\alpha_{3}\rangle$ is a nonzero multiple of the jumping generator $\epsilon$.
\end{proof}

\subsection{The integral homology $H^{*}(\Symp(\tilde{M}^{0}_{\mu,c});\Z)$ when $\mu\in(1,2)$ and $\lambda\leq c$}

We now use the description of the stratification of the space of tamed almost complex structures $\tilde{\jj}_{\mu,c}=\jj(\tMuco)$ given in Section~\ref{SectionStratification} to prove that $H^{*}(\Symp(\tilde{M}^{0}_{\mu,c});\Z)$ has no torsion.

For that range of parameters, the contractible space $\tilde{\jj}_{\mu,c}=\tilde{\jj}_{01}$ is the union of the open, dense, and connected stratum $\tilde{\jj}_{0}$ with the codimension $2$, co-oriented, closed submanifold $\tilde{\jj}_{1}$. Hence, a tubular neighborhood $\Nn$ of $\tilde{\jj}_{1}$ is homeomorphic to an oriented disc bundle $D^{2}\to\Nn\to\tilde{\jj}_{1}$. Let $\F$ denote either $\Q$ or a finite field $\F_{p}$ for $p\geq 2$ prime.

\begin{lemma}\label{Isomorphisme2Strates}
There is an isomorphism
$H_{i+1}(\tilde{\jj}_{0};\F)\simeq H_i(\tilde{\jj}_{1};\F)$, $\forall\, i\geq 0$.
\end{lemma}

\begin{proof}
Since the space $\tilde{\jj}_{01}$ is contractible, the long exact sequence of the pair  $(\tilde{\jj}_{01},\tilde{\jj}_{0})$ shows that the groups $H_{i+2}(\tilde{\jj}_{01},\tilde{\jj}_{0})$ and $H_{i+1}(\tilde{\jj}_{0})$ are isomorphic. Let write $\xi$ for the intersection $\tilde{\jj}_{0}\cap\Nn$. Excising $\tilde{\jj}_{01}-\Nn$ from $(\tilde{\jj}_{01},\tilde{\jj}_{0})$, we see that $H_{i+2}(\tilde{\jj}_{01},\tilde{\jj}_{0})\simeq H_{i+2}(\Nn, \xi)$. By the Thom isomorphism theorem, there is an isomorphism
$H_{i+2}(\Nn), \xi) \simeq H_{i}(\tilde{\jj}_{1})$ for all $i\geq 0$.
\end{proof}

Each stratum $\tilde{\jj}_{i}$ is homotopy equivalent to the quotient $\Symp(\tMuco)/\tilde{T}_{i}$ and we have two homotopy fibrations $\tilde{T}_{i} \to \Symp(\tMuco)\to \tilde{\jj}_{i}$.

\begin{prop}\label{InclusionEnHomologie}
The inclusion $\tilde{T^2_i}\hookrightarrow \Symp(\tMuco)$ induces an injection in homology for all $i \geq 0$ and $\mu \geq 1$, and for any coefficient field $\F$.
\end{prop}

\begin{proof}
First note that the inclusion $\tilde{T}^2_i\hookrightarrow \Symp(\tMuco)$ induces a morphism of Hopf algebras. Since $H_{*}(\tilde{T}^2_i;\F)$ is primitively generated and commutative, it is enough to show that its submodule of primitive elements, which is simply $H_1(\tilde{T}^2_i;\F)$, is mapped injectively into $H_1(\Symp(\tMuco);\F)$.

The same proof as in Lemma~\ref{RetractionSymp} shows that the space $\Symp(\tMuco)$ is homotopy equivalent to $\Symp(\tMuco, E \cup (F-E))$. One therefore gets a map from $\Symp(\tMuco) \to \Symp(\tMuco, E \cup (F-E)) \to \Ee(W_{i}, F_i, *)$, the space of homotopy self-equivalences of $W_i$ that preserve a fiber $F_i$, and $*$, the point of intersection between $F_i$ and the section of self-intersection $-i$. Here the last map is the one induced by blowing down an exceptional curve in class $E$ when $i=1$ and $F-E$ when $i=2$. Consider the maps
$$
\Ee(W_{i}, F_i, *) \to \Ee(F_i, *)
$$
given by the restriction to the fiber $F_i$ and 
$$
\Ee(W_{i}, F_i, *) \to \Ee(S^2, *)
$$
given by first embedding $S^{2}$ in $W_{i}$ as the zero section, applying a homotopy equivalence, and then projecting back to the base of the Hirzebruch fibration. These two maps show that the images in $\pi_{1}(\Symp(\tMuco))$ of the two generators of $\pi_{1}(\tilde{T}_{i})$ have infinite order and are independent. Since $\Symp(\tMuco)$ is a topological group, $\pi_{1}(\Symp(\tMuco))=H_{1}(\Symp(\tMuco);\Z)$ and the lemma follows readily.
\end{proof}

\begin{cor}\label{ScindementHomologique}
There is an isomorphism of graded vector spaces $H_*(\Symp(\tilde{M}_{\mu,c});\F) \simeq H_*(\tilde{\jj}_i;\F)\otimes H_*(\tilde{T}^2_i;\F)$ and an isomorphism of graded algebras $H^*(\Symp(\tilde{M}_{\mu,c});\F) \simeq H^*(\tilde{\jj}_i;\F)\otimes H^*(\tilde{T}^2_i;\F)$. Consequently, $H^*(\tilde{\jj}_0;\F)\simeq H^*(\tilde{\jj}_1;\F)$.
\end{cor}

\begin{proof}
This follows from the Leray-Hirsh theorem applied to the homotopy fibrations
\[\tilde{T}_{i} \to \Symp(\tMuco)\to \tilde{\jj}_{i}\] 
and from the fact that $H^{*}(\tilde{T}^2_i;\F)$ is a free algebra.
\end{proof}

\begin{prop} \label{IntegralCohomology} 
When $1<\mu<2$ and $\lambda\leq c$, $H_{*}(\Symp(\tMuco);\Z)$ has no torsion.
\end{prop}
\begin{proof}
Combining Lemma~\ref{Isomorphisme2Strates} with Corollary~\ref{ScindementHomologique} and using the fact that $\tilde{\jj}_i$ is connected, we get $H^q(\tilde{\jj}_i;\F)\simeq\F$, for all $q\geq 0$. Putting this back into the isomorphism
\[
H_{*}(\Symp(\tilde{M}_{\mu,c});\F) \simeq H_{*}(\tilde{\jj}_i;\F)\otimes H_{*}(\tilde{T}^2_i;\F)
\] 
this shows that the homology $H_{*}(\Symp(\tilde{M}_{\mu,c});\F)$ is the tensor product $H_{*}(\Symp(\tilde{M}_{\mu,c});\Z)_{\text{free}}\otimes \F$ and this implies that the integral homology (and cohomology) has no torsion.
\end{proof}
From Proposition~\ref{IsomorphismFundamentalGroups}, we get
\begin{cor}\label{CorPi1}
For all $\mu>1$ and $c\in(0,1)$, $\pi_{1}(\Symp(\tMuco))\simeq \Z\oplus\Z\oplus\Z$.
\end{cor}

\begin{remark}\label{RemarkHomotopyEquivalence}
When $\mu\in(1,2)$ and $\lambda\leq c$, it can be shown, using arguments similar to the ones used in Anjos~\cite{An} and Anjos-Granja~\cite{AG}, that $\Symp(\tMuco)$ is homotopy equivalent (as a space) to $\Omega S^{3}\times S^{1}\times S^{1}\times S^{1}$. 
\end{remark}

\section{Spaces of embeddings}\label{SectionEmbedddings}
We now investigate the homotopy type of the space of symplectic embeddings $\IEmb(B_{c},\Mui)$.

\begin{proof}[Proof of Theorem~\ref{HomotopyOfEmbeddings}]
Recall that McDuff~\cite{MD-Isotopy} showed that the embedding spaces $\IEmb(B_{c},\Mui)$ are always connected. Consequently, we have a homotopy fibration
\begin{equation}\label{FibrationPrincipale}
\Symp(\tMuci,\Sigma)\to \Symp(\Mui)\to \IEmb(B_c,\Mui)
\end{equation}
When $c<\lambda$, Theorem~\ref{HomotopyOfSymp} states that the fiber $\Symp(\tMuci,\Sigma)$ is homotopy equivalent to the symplectic stabilizer of a point in $\Mui$. Consequently, $\IEmb(B_c,\Mui)$ has the same homotopy type as $\Mui$.

When $\lambda\leq c<\mu$, the rational homotopy groups of $\Symp(\tMuci,\Sigma)$ and $\Symp(\Mui)$ are given in Corollary~\ref{RationalHomotopy} and Theorem~\ref{AbreuMcDuff}. The only ambiguous part in the rational homotopy sequence of the above fibration is
\[
\to 0 \to\pi_{2}(\IEmb)\to\pi_{1}(\Symp(\tMuco))\to\pi_{1}(\Symp(\Muo))\to\pi_{1}(\IEmb)\to 0
\]
But since the blow-up of the generator of the group $\pi_{1}(\Symp(\Muo))\simeq \Q$ is one of the generators of $\pi_{1}(\Symp(\tMuco))\simeq \Q^{3}$, the map $\pi_{1}(\Symp(\tMuco))\to\pi_{1}(\Symp(\Muo))$ is surjective and, therefore, $\pi_{1}(\IEmb)=0$ and $\pi_{2}(\IEmb)\simeq \Q$. The rest of the computation of $\pi_{*}(\IEmb)$ is straightforward.
\end{proof}

\begin{proof}[Proof of Theorem~\ref{HomotopyOfEmbeddingsSpecialCase}]
This follows directly from Theorem~\ref{HomotopyOfSympSpecialCase} applied to the fibration~(\ref{FibrationPrincipale}).
\end{proof}

\subsection{Symplectic balls in $\CP^{2}$}

We now apply the relative framework of Section~\ref{SectionHomotopySpecialCase} to determine the homotopy type of the space of one or two disjoint symplectic balls in $\CP^{2}$. Recall that if we normalize $\CP^{2}$ so that the area of a line $L$ is $1$, then the embedding spaces $\Emb(B_{c},\CP^{2})$ and $\Emb(B_{c_{1}}\sqcup B_{c_{2}},\CP^{2})$ are nonempty if, and only if, $\delta\in (0,1)$ and $\delta_{1}+\delta_{2}\in(0,1)$, see~\cite{MP}. By the work of McDuff~\cite{MD-Isotopy}, we know that these spaces are connected.

\begin{proof}[Proof of Theorem~\ref{EmbeddingsInCP2}]
We first consider the embedding of one ball in $\CP^{2}$. Let $\Sigma$ denote the exceptional divisor in $\tCP^{2}$ and consider the fibrations
\begin{equation}\label{FibrationCP2}
\Symp(\tCP^{2},\Sigma)\to \Symp(\CP^{2})\to \IEmb(B_\delta,\CP^{2})
\end{equation}
and
\begin{equation}\label{SecondFibrationCP2}
\Symp(\tCP^{2},\Sigma)\to \Symp(\tCP^{2})\to \Cc([\Sigma])
\end{equation}
where $\Cc([\Sigma])$ is the space of embedded symplectic spheres in the class of the exceptional divisor. The blow-up manifold $\tCP^{2}$ is conformally symplectomorphic to the nontrivial bundle $\Mul$ with $\mu=\delta/(1-\delta)>0$. The space $\jj_{\mu}^{1}$ of all tame almost complex structures on $\Mul$ has a stratification induced by the degeneration type of the class $[\Sigma]=B^{1}$ (which, for $\mu>k\in \N$, can be represented by cusp curves containing a component in the class $B^{1}-kF^{1}$). In particular, the open stratum $\jj_{\mu,0}$ (which is made of almost complex structures $J$ for which the exceptional class is represented by a $J$-holomorphic embedded sphere) is homotopy equivalent to the space $\Cc([\Sigma])$. By the work of Abreu and McDuff~\cite{AM}, we know that the open stratum is also homotopy equivalent to $\Symp(\Mul)/\U(2)$, where the elements of $\U(2)$ are lifts to $\Mul$ of complex automorphisms of $\CP^{2}$ fixing a point. It follows that, for any $\mu=\delta/(1-\delta)$, the stabilizer $\Symp(\tCP^{2},\Sigma)$ of $\Sigma$ in the fibration (\ref{SecondFibrationCP2}) is homotopy equivalent to $\U(2)$. Since $
\Symp(\CP^{2})\simeq \PU(3)$, this implies that the fibration~(\ref{FibrationCP2}) is in turn homotopy equivalent to the standard fibration
\[\U(2)\to\PU(3)\to\CP^{2}\]
and this conludes the proof of the first statement in Theorem~\ref{EmbeddingsInCP2}.

We now consider the embedding of two disjoint balls $\Bb:=B_{\delta_{1}}\sqcup B_{\delta_{2}}$ in $\CP^{2}$, where we suppose that $\delta_{1}\geq\delta_{2}$. Let denote by $\Dd=\Sigma_{1}\sqcup\Sigma_{2}$ the union of the exceptional divisors in the $2$-fold blow-up $X_{\delta_{1},\delta_{2}}$ and consider the fibrations
\begin{equation}\label{Fibration2BallsCP2}
\Symp(X_{\delta_{1},\delta_{2}},\Dd)\to \Symp(\CP^{2})\to \IEmb(\Bb,\CP^{2})
\end{equation}
\begin{equation}\label{RelativeFibration2Balls}
\Symp(X_{\delta_{1},\delta_{2}},\Dd)\cap\Diff_{0}(X_{\delta_{1},\delta_{2}},\Dd)\to \Diff_{0}(X_{\delta_{1},\delta_{2}},\Dd)\to \Ss_{0}(\om,\Dd)\simeq\Aa(\om,\Dd)
\end{equation}
This time, $X_{\delta_{1},\delta_{2}}$ is conformally symplectomorphic to $\tMuco$ with $\mu=(1-\delta_{2})/(1-\delta_{1}) \geq 1$ and $c=(1-\delta_{1}-\delta_{2})/(1-\delta_{1})\in(0,1)$ via a diffeomorphism which sends the exceptional classes $[L]-[\Sigma_{1}]-[\Sigma_{2}]$, $[\Sigma_{1}]$, and $[\Sigma_{2}]$ to, respectively, $E$, $B-E$, and $F-E$. In particular, $\Dd$ is a configuration of disjoint exceptional spheres in classes $B-E$ and $F-E$. Let write $\Cc(\Dd)$ for the set of all such configurations. Because, for any $J\in\jjuco$, the class $F-E$ is represented by a unique embedded $J$-holomorphic sphere, $\Cc(\Dd)$ is homotopy equivalent to the space $\Cc(B-E)$ of all exceptional symplectic spheres   in class $B-E$. As explained before, $\Cc(B-E)$ is homotopy equivalent to the open stratum $\jj_{\mu,c,0}$ which is itself homotopy equivalent to the homogeneous space $\Symp(\tMuco)/T_{0}$. It follows that $\Symp(X_{\delta_{1},\delta_{2}},\Dd)\simeq \Symp(\tMuco,\Dd)$ is homotopy equivalent to the torus $T_{0}$. In particular, $\Symp(X_{\delta_{1},\delta_{2}},\Dd)$ is connected.

Because $\Aa(\om,\Dd)$ is contained in the open stratum of $\Aa_{\mu,c}$, it is independent of $\mu$ and $c$. As before, it follows that the homotopy type of the group $\Symp(X_{\delta_{1},\delta_{2}},\Dd)$ is independent of $\delta_{1}$ and $\delta_{2}$. Considering the limit as $\delta_{1}$ and $\delta_{2}=1-\delta_{1}-\epsilon$ go to zero shows that $\Symp(X_{\delta_{1},\delta_{2}},\Dd)$ is homotopy equivalent to the stabilizer of two points in $\CP^{2}$. This implies that $\IEmb(\Bb,\CP^{2})$ is homotopy equivalent to the space of ordered configurations of two points in $\CP^{2}$. 
\end{proof}


\begin{ackn}
I am pleased to thank F. Lalonde, D. McDuff, and S. Anjos for helpful discussions, and also an anonymous referee for pointing out an error in my original treatment of Theorem~\ref{StabilityTheorem} and for suggesting improvements to the exposition. Part of this work was carried out while I was a postdoctoral
fellow at the Department of Mathematics of the University of Toronto and I
take this opportunity to thank the staff and faculty for providing
a stimulating environment. I also thank the Fields Institute where this work was completed.
\end{ackn}


\end{document}